\documentstyle{article}

\def\eqk{\equiv_{\cal K}}

\newcommand{\df}{\mbox{\scriptsize{\it df}}}

\def\dk#1{\mbox{$\lfloor {#1} \rfloor$}}

\def\gk#1{\mbox{$\lceil {#1} \rceil$}}

\def\eqk{\equiv_{\cal K}}

\def\Lo{{\mbox{${\cal L}_\omega$}}}

\def\Ko{{\mbox{${\cal K}_\omega$}}}

\def\Mat{{\mbox{${\mbox{\bf Mat}}_{\cal F}$}}}

\newcommand{\mj}{\mbox{\bf 1}}

\newcommand{\str}{\rightarrow}

\newcommand{\N}{{\mbox{\bf N}}}

\newcommand{\Z}{{\mbox{\bf Z}}}

\def\SKo{{\mbox{${\cal SK}_\omega$}}}

\def\SKn{{\mbox{${\cal SK}_n$}}}

\def\SK{{\mbox{$\cal SK$}}}

\def\SJo{{\mbox{${\cal SJ}_\omega$}}}

\def\SJ{{\mbox{$\cal SJ$}}}

\def\cirk{\,{\raisebox{.3ex}{\tiny $\circ$}}\,}

\newcommand{\De}{\mbox{$\cal D$}}

\def\prop#1#2{\vspace{2ex} \noindent{\sc #1.} {\it #2} \par \vspace{2ex}}

\def\dkz{\noindent{\emph{Proof.  }}}

\def\qed{\hfill $\dashv$}

\def\pl{\!+\!}

\def\mn{\!-\!}

\def\nav#1#2{\parshape=2 1em 32em 4em 29em \smallskip
\noindent{\makebox[3em][l]{#1}}{#2}\par\smallskip}

\def\So{{\mbox{${\cal S}_\omega$}}}

\begin{document}

\title
{\bf Symmetric Self-Adjunctions \\and Matrices}

\author{
{\sc Kosta Do\v sen} and {\sc Zoran Petri\' c}
\\[.3cm]{\small Mathematical Institute, SANU}
\\[-.1cm]{\small Knez Mihailova 35, P.O. Box 367}
\\[-.1cm]{\small 11001 Belgrade, Serbia}
\\[-.1cm]{\small email: \{kosta, zpetric\}@mi.sanu.ac.rs}}

\date{ }
\maketitle

\begin{abstract}
\noindent It is shown that the multiplicative monoids of Brauer's
centralizer algebras generated out of the basis are isomorphic to
monoids of endomorphisms in categories where an endofunctor is
adjoint to itself, and where, moreover, a kind of symmetry
involving the self-adjoint functor is satisfied. As in a previous
paper, of which this is a companion, it is shown that such a
symmetric self-adjunction is found in a category whose arrows are
matrices, and the functor adjoint to itself is based on the
Kronecker product of matrices.

\vspace{0.3cm}

\noindent{\it Mathematics Subject Classification} ({\it 2000}):
18A40, 57M99, 20B30

\vspace{0.1cm}

\noindent{\it Keywords}: adjunction, matrix representation,
symmetric groups, Brauer's centralizer algebras

\end{abstract}

\section{Introduction}

This paper is a companion to \cite{DP03a}, where it was shown that
the category \Mat\ whose objects are dimensions of
finite-dimensional vector spaces and whose arrows are matrices is
an example of an adjoint situation where a functor is adjoint to
itself. In \Mat\ this functor is based on the Kronecker product of
matrices. This self-adjunction underlies the orthogonal group case
of Brauer's representation of Brauer's algebras, which can be
restricted to the Temperley-Lieb subalgebras of Brauer's algebras
(see \cite{B37}, \cite{W88}, \cite{J94}).

In this paper we show that this self-adjunction of \Mat\ is
specific in that it involves a kind of symmetry (see Section~8).
This symmetry is the reason why the self-adjunction of \Mat\ is of
the $\cal K$ kind and not of the $\cal L$ kind, in the terminology
of \cite{DP03a} (see the end of Section~8; cf.\ Section 2).
Roughly speaking, in the diagrammatic representation, the $\cal L$
kind takes account of the regions where the circles occur, while
the $\cal K$ kind just takes account of the number of circles, and
does not take account of where they occur. In the $\cal J$ kind,
one does not take account of circles at all.

The first part of the paper (Sections 2-7) is devoted to
axiomatizing (presenting by generators and relations) monoids for
which we will show later that they are engendered by categories
involved in symmetric self-adjoint situations. These monoids are
closely related to multiplicative monoids of Brauer's algebras. We
establish that our monoids are isomorphic to monoids of diagrams
which differ from the tangles of knot theory by having crossings
in which \emph{under} and \emph{over} are disregarded. These
crossings correspond to the symmetry of self-adjunction mentioned
above. In the first part of the paper (Section~7) we also
establish a property of our monoids that we call
\emph{maximality}, on which we rely later for the interpretation
in \Mat. (There is an alternative to this approach based on
\cite{DP03c} and \cite{DP03d}.)

In the second part (Sections 8-10) we introduce formally and
gradually our notion of symmetric self-adjunction, starting from
more general notions---ultimately from the notion of adjunction.
We favour this gradual approach in order to make clearer the
articulation of this notion. In this part we also connect freely
generated symmetric self-adjunctions with the monoids of the first
part.

In the third part (Sections 11-13) we survey a number of notions
of category related to symmetric self-adjunction, but more
general. Such are the notions of symmetric monoidal closed
category and compact closed category. The notion that corresponds
exactly to symmetric self-adjunction is a notion we call
\emph{subsided} category. It is shown that the category of
symmetric self-adjunction freely generated by a single object is
isomorphic to the subsided category $\cal S$ freely generated by a
single object. In $\cal S$ the objects may be conceived as
functors and the tensor product as composition of functors. If in
the generating set of $\cal S$ we have more than one object, then
we fall upon a generalization of our notion of symmetric
self-adjunction, which we call \emph{multiple} symmetric
self-adjunction (see Section 13).

The final part of the paper (Section 14) is about the symmetric
self-adjunction of \Mat. In \Mat\ we have as subcategories
isomorphic copies of freely generated symmetric self-adjunctions.

\section{The monoids \SKo\ and \SJo}

The monoid \SKo\ has for every ${k\in \N^+=\N\!-\!\{0\}}$ a
generator \dk{k}, called a \emph{cup term}, a generator \gk{k},
called a \emph{cap term}, and a generator $\sigma_k$, called a
\emph{crossing term}. (We prefer the name \emph{crossing} to the
name \emph{transposition} because of the pictorial interpretation
of $\sigma_k$ in the next section.) The \emph{terms} of \SKo\ are
defined inductively by stipulating that the generators and \mj\
are terms, and that if $t$ and $u$ are terms, then $(tu)$ is a
term. As usual, we will omit the outermost parentheses of terms,
and in the presence of associativity we will omit all parentheses,
since they can be restored as we please.

The monoid \SKo\ is freely generated from the generators above so
that the following equations hold between terms of \SKo\ for
$j\leq k$:

\begin{tabbing}
\hspace{1em}\emph{monoid equations}:
\\*[1ex]
\hspace{11em}\= (1)\quad\= $\mj t=t\mj=t$,
\\*[1ex]
\> (2)\> $t(uv)=(tu)v$;
\\[2ex]
\hspace{1em}\emph{cup-cap equations}:
\\*[1ex]
\hspace{2em}\= (\emph{cup})\hspace{3.5em}\= $\dk{k}\dk{j} =
\dk{j}\dk{k+2}$,\hspace{2.5em}\= (\emph{cap})\hspace{3.5em}\=
$\gk{j}\gk{k} = \gk{k+2}\gk{j}$,
\\*[1ex]
\> (\emph{cup-cap}~1)\> $\dk{k+2}\gk{j} = \gk{j}\dk{k}$, \>
(\emph{cap-cup}~1)\> $\dk{j}\gk{k+2} = \gk{k}\dk{j}$,
\\*[1ex]
\hspace{11em}(\emph{cup-cap})\quad $\dk{i}\gk{i+1} =\mj$;
\\[2ex]
\hspace{1em}$\sigma$ \emph{equations}:
\\*[1ex]
\hspace{11em}\= ($\sigma$)\quad\=
$\sigma_{k+2}\sigma_j=\sigma_j\sigma_{k+2}$,
\\*[1ex]
\> ($\sigma$2)\> $\sigma_i\sigma_i=\mj$,
\\*[1ex]
\> ($\sigma$3)\>
$\sigma_{i+1}\sigma_i\sigma_{i+1}=\sigma_i\sigma_{i+1}\sigma_i$;
\\[2ex]
\hspace{1em}$\sigma$-\emph{cup equations}:\hspace{10em}
$\sigma$-\emph{cap equations}:
\\*[1ex]
\hspace{2em}\= (\emph{cup})\hspace{3em}\= $\dk{k}\dk{j} =
\dk{j}\dk{k+2}$,\hspace{3em}\= (\emph{cap})\hspace{3em}\=
$\gk{j}\gk{k} = \gk{k+2}\gk{j}$,\kill \> ($\sigma$-\emph{cup}~1)\>
$\dk{k+2}\sigma_j=\sigma_j\dk{k+2}$, \> ($\sigma$-\emph{cap}~1)\>
$\sigma_j\gk{k+2}=\gk{k+2}\sigma_j$,
\\*[1ex]
\> ($\sigma$-\emph{cup}~2)\> $\dk{j}\sigma_{k+2}=\sigma_k\dk{j}$,
\> ($\sigma$-\emph{cap}~2)\> $\sigma_{k+2}\gk{j}=\gk{j}\sigma_k$,
\\*[1ex]
\> ($\sigma$-\emph{cup}~3)\> $\dk{i}\sigma_i=\dk{i}$, \>
($\sigma$-\emph{cap}~3)\> $\sigma_i\gk{i}=\gk{i}$,
\\*[1ex]
\> ($\sigma$-\emph{cup}~4)\>
$\dk{i+1}\sigma_i=\dk{i}\sigma_{i+1}$, \>
($\sigma$-\emph{cap}~4)\> $\sigma_i\gk{i+1}=\sigma_{i+1}\gk{i}$.
\end{tabbing}

To understand these equations it helps to have in mind their
diagrammatic interpretation of the next section. The monoid and
cup-cap equations are the equations of the monoid \Lo\ of
\cite{DP03a} and \cite{DP03b}, while the monoid and $\sigma$
equations are the well known equations of symmetric groups (see
\cite{CM57}, Chapter~6.2, p.~64).

The equation

\[
\dk{k}\gk{k}=\dk{k+1}\gk{k+1}
\]

\noindent holds in \SKo\ since we have

\begin{tabbing}
\hspace{6em}$\dk{k}\gk{k}$ \=
$=\dk{k}\sigma_{k+1}\sigma_{k+1}\gk{k}$,\quad \= with ($\sigma$2)
\\[1ex]
\> $=\dk{k+1}\gk{k+1}$,\> with ($\sigma$-\emph{cup}~4),
($\sigma$-\emph{cap}~4) and ($\sigma$2).
\end{tabbing}

\noindent The monoid \Ko\ of \cite{DP03a} and \cite{DP03b}, in
which $\sigma_i$ is lacking, is a submonoid of \SKo.

Let us call terms of the form $\dk{k}\gk{k}$ \emph{circles}. Then,
by the equation we have just derived, we have
$\dk{k}\gk{k}=\dk{l}\gk{l}$, which means that we have only one
circle, which we designate by $c$.

We can easily derive the equation

\[
(c)\quad tc=ct
\]

\noindent for $t$ being $\dk{k}$, $\gk{k}$ or $\sigma_i$, which
yields $(c)$ for any term $t$.

The monoid \SJo\ is defined as \SKo\ save that we have the
additional equation

\[
(c1)\quad c=\mj.
\]

Let $c^0$ be the empty sequence, and let $c^{n+1}$ be ${c^nc}$.
Then it is easy to show by induction on the length of derivation
that if $t=u$ in \SJo, then for some $n,m\in \N$ we have
$c^nt=c^mu$ in \SKo. (The converse implication holds trivially.)
Later (after Section~6) we will be able to establish that if
${t=u}$ holds in \SJo\ but not in \SKo, then we have ${c^nt=c^mu}$
in \SKo\ for $n\neq m$ (cf.\ \cite{DP03a}, end of Section 12).

\section{\SK\ and \SJ\ diagrams}

A \emph{coupling diagram} is a partition of ${\Z\!-\!\{0\}}$ into
two-element subsets. (Alternatively, we could define a coupling
diagram by an involution of ${\Z\!-\!\{0\}}$ without fixed-points;
a coupling diagram may be conceived as a \emph{split equivalence}
in the sense of \cite{DP03c} and \cite{DP03d}.) The two-element
sets of integers that belong to a coupling diagram are called
\emph{threads} of the diagrams. When in a thread ${\{a,b\}}$ the
members $a$ and $b$ have a different sign, the thread is called
\emph{transversal}. A thread ${\{a,b\}}$, which is not
transversal, is a \emph{cup} when $a$ and $b$ are both positive,
and it is a \emph{cap} when they are both negative.

We denote threads systematically as follows. A transversal thread
is denoted by ${\{n,-m\}}$ for ${n,m\in \N^+}$, a cup is denoted
by ${\{n,m\}}$ with the assumption that ${0<n<m}$, and a cap is
denoted by ${\{-n,-m\}}$, again with the assumption that
${0<n<m}$.

Consider coupling diagrams such that for some ${n,m\in \N^+}$ for
every ${k\in \N}$ we have a transversal thread ${\{n\pl k,-(m\pl
k)\}}$. Such coupling diagrams will be called \SJ\
\emph{diagrams}, and ${(n,m)}$ is a \emph{type} of this diagram.
The type of an \SJ\ diagram is not unique, since an \SJ\ diagram
of type ${(n,m)}$ is also of type ${(n\pl k,m\pl k)}$. It is clear
that an \SJ\ diagram has finitely many cups and caps. A coupling
diagram without cups and caps corresponds to a permutation of
$\N^+$, and an \SJ\ diagram without cups and caps corresponds to a
permutation of an initial segment of $\N^+$. In a type ${(n,m)}$
of such an \SJ\ diagram we must have ${n=m}$.

An \SK\ \emph{diagram} is a pair ${(D,n)}$ where $D$ is an \SJ\
diagram and ${n\in \N}$ (we imagine that $n$ is the number of
\emph{circular components} of the diagram; see the picture of
composition below). We draw an \SJ\ diagram by putting positive
integers on a line, called the \emph{top line}, and negative
integers in reverse order on a line below, called the \emph{bottom
line}, and by connecting for each thread ${\{a,b\}}$, the numbers
$a$ and $b$ by a, not necessarily straight, line segment in
between the top and bottom lines. For example, the \SJ\ diagram

\begin{tabbing}
$\{\{1,-3\},\{2,3\},\{4,-4\},\{5,-1\},\{-2,-6\},\{-5,-7\},\{6,-8\},
\ldots,$
\\
\`$\{6\pl k,-(8\pl k)\},\ldots\}$
\end{tabbing}

\noindent is drawn as follows:

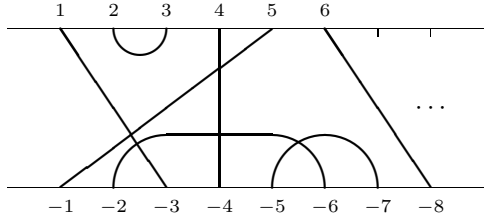
\begin{center}
\begin{picture}(180,100)
{\linethickness{0.05pt} \put(0,20){\line(1,0){180}}
\put(0,80){\line(1,0){180}} \put(140,80){\line(0,-1){3}}
\put(160,80){\line(0,-1){3}}}

\put(20,85){\makebox(0,0)[b]{\scriptsize $1$}}
\put(40,85){\makebox(0,0)[b]{\scriptsize $2$}}
\put(60,85){\makebox(0,0)[b]{\scriptsize $3$}}
\put(80,85){\makebox(0,0)[b]{\scriptsize $4$}}
\put(100,85){\makebox(0,0)[b]{\scriptsize $5$}}
\put(120,85){\makebox(0,0)[b]{\scriptsize $6$}}

\put(20,15){\makebox(0,0)[t]{\scriptsize $-1$}}
\put(40,15){\makebox(0,0)[t]{\scriptsize $-2$}}
\put(60,15){\makebox(0,0)[t]{\scriptsize $-3$}}
\put(80,15){\makebox(0,0)[t]{\scriptsize $-4$}}
\put(100,15){\makebox(0,0)[t]{\scriptsize $-5$}}
\put(120,15){\makebox(0,0)[t]{\scriptsize $-6$}}
\put(140,15){\makebox(0,0)[t]{\scriptsize $-7$}}
\put(160,15){\makebox(0,0)[t]{\scriptsize $-8$}}

\thicklines \put(20,20){\line(4,3){80}}
\put(60,20){\line(-2,3){40}} \put(80,20){\line(0,1){60}}
\put(160,20){\line(-2,3){40}}

\put(160,50){\makebox(0,0){$\cdots$}}

\put(50,80){\oval(20,20)[b]} \put(80,20){\oval(80,40)[t]}
\put(120,20){\oval(40,40)[t]}
\end{picture}
\end{center}

\noindent This \SJ\ diagram is of type ${(6,8)}$, ${(7,9)}$, etc.

Let the \emph{unit} \SJ\ diagram be the coupling diagram
${\{\{n,-n\}\mid n\in \N^+\}}$, which is drawn as follows:

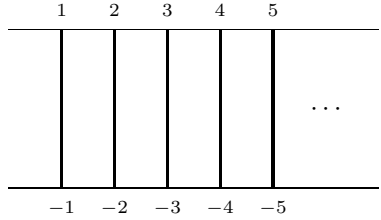
\begin{center}
\begin{picture}(140,100)
{\linethickness{0.05pt} \put(0,20){\line(1,0){140}}
\put(0,80){\line(1,0){140}}}

\put(20,85){\makebox(0,0)[b]{\scriptsize $1$}}
\put(40,85){\makebox(0,0)[b]{\scriptsize $2$}}
\put(60,85){\makebox(0,0)[b]{\scriptsize $3$}}
\put(80,85){\makebox(0,0)[b]{\scriptsize $4$}}
\put(100,85){\makebox(0,0)[b]{\scriptsize $5$}}

\put(20,15){\makebox(0,0)[t]{\scriptsize $-1$}}
\put(40,15){\makebox(0,0)[t]{\scriptsize $-2$}}
\put(60,15){\makebox(0,0)[t]{\scriptsize $-3$}}
\put(80,15){\makebox(0,0)[t]{\scriptsize $-4$}}
\put(100,15){\makebox(0,0)[t]{\scriptsize $-5$}}

\thicklines \put(20,20){\line(0,1){60}}
\put(40,20){\line(0,1){60}} \put(60,20){\line(0,1){60}}
\put(80,20){\line(0,1){60}} \put(100,20){\line(0,1){60}}

\put(120,50){\makebox(0,0){$\cdots$}}
\end{picture}
\end{center}

\noindent and is denoted by $I$. The unit \SK\ diagram is
${(I,0)}$.

The \SJ\ diagram that is the composition $D_2\cirk D_1$ of two
\SJ\ diagrams is defined pictorially by identifying the bottom
line of $D_1$ with the top line of $D_2$, so that $-n$ and $n$ are
the same point. The top line of $D_2\cirk D_1$ is the top line of
$D_1$, and the bottom line of $D_2\cirk D_1$ is the bottom line of
$D_2$. For example, for $D_1$ being the \SJ\ diagram taken as an
example above, and $D_2$ the \SJ\ diagram

\[
\{\{1,-1\},\{2,7\},\{3,-3\},\{4,-2\},\{5,6\},\{8,-4\},
\ldots,\{8\pl n,-(4\pl n)\},\ldots\},
\]

\noindent the composition $D_2\cirk D_1$ is

\[
\{\{1,-3\},\{2,3\},\{4,-2\},\{5,-1\},\{6,-4\}, \ldots,\{6\pl
n,-(4\pl n)\},\ldots\},
\]

\noindent which is clear from the following pictures:

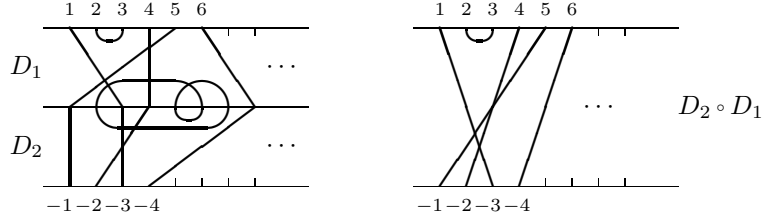
\begin{center}
\begin{picture}(240,100)
{\linethickness{0.05pt} \put(0,20){\line(1,0){100}}
\put(0,50){\line(1,0){100}} \put(0,80){\line(1,0){100}}
\put(140,20){\line(1,0){100}} \put(140,80){\line(1,0){100}}
\put(70,80){\line(0,-1){3}} \put(80,80){\line(0,-1){3}}
\put(210,80){\line(0,-1){3}} \put(220,80){\line(0,-1){3}}
\put(50,20){\line(0,1){3}} \put(60,20){\line(0,1){3}}
\put(70,20){\line(0,1){3}} \put(80,20){\line(0,1){3}}
\put(190,20){\line(0,1){3}} \put(200,20){\line(0,1){3}}
\put(210,20){\line(0,1){3}} \put(220,20){\line(0,1){3}}}

\put(10,85){\makebox(0,0)[b]{\scriptsize $1$}}
\put(20,85){\makebox(0,0)[b]{\scriptsize $2$}}
\put(30,85){\makebox(0,0)[b]{\scriptsize $3$}}
\put(40,85){\makebox(0,0)[b]{\scriptsize $4$}}
\put(50,85){\makebox(0,0)[b]{\scriptsize $5$}}
\put(60,85){\makebox(0,0)[b]{\scriptsize $6$}}

\put(6,15){\makebox(0,0)[t]{\scriptsize $-1$}}
\put(17,15){\makebox(0,0)[t]{\scriptsize $-2$}}
\put(28,15){\makebox(0,0)[t]{\scriptsize $-3$}}
\put(39,15){\makebox(0,0)[t]{\scriptsize $-4$}}

\put(150,85){\makebox(0,0)[b]{\scriptsize $1$}}
\put(160,85){\makebox(0,0)[b]{\scriptsize $2$}}
\put(170,85){\makebox(0,0)[b]{\scriptsize $3$}}
\put(180,85){\makebox(0,0)[b]{\scriptsize $4$}}
\put(190,85){\makebox(0,0)[b]{\scriptsize $5$}}
\put(200,85){\makebox(0,0)[b]{\scriptsize $6$}}

\put(146,15){\makebox(0,0)[t]{\scriptsize $-1$}}
\put(157,15){\makebox(0,0)[t]{\scriptsize $-2$}}
\put(168,15){\makebox(0,0)[t]{\scriptsize $-3$}}
\put(179,15){\makebox(0,0)[t]{\scriptsize $-4$}}

\thicklines \put(10,50){\line(4,3){40}}
\put(30,50){\line(-2,3){20}} \put(40,50){\line(0,1){30}}
\put(80,50){\line(-2,3){20}}

\put(10,20){\line(0,1){30}} \put(30,20){\line(0,1){30}}
\put(20,20){\line(2,3){20}} \put(40,20){\line(4,3){40}}

\put(150,20){\line(2,3){40}} \put(160,20){\line(1,3){20}}
\put(170,20){\line(-1,3){20}} \put(180,20){\line(1,3){20}}

\put(90,65){\makebox(0,0){$\cdots$}}
\put(90,35){\makebox(0,0){$\cdots$}}

\put(210,50){\makebox(0,0){$\cdots$}}

\put(25,80){\oval(10,10)[b]} \put(40,50){\oval(40,20)[t]}
\put(60,50){\oval(20,20)[t]}

\put(45,50){\oval(50,16)[b]} \put(55,50){\oval(10,10)[b]}

\put(165,80){\oval(10,10)[b]}

\put(0,65){\makebox(0,0)[r]{$D_1$}}
\put(0,35){\makebox(0,0)[r]{$D_2$}}

\put(240,50){\makebox(0,0)[l]{$D_2\cirk D_1$}}

\end{picture}
\end{center}

We have in the picture on the left a closed loop made of the two
caps of $D_1$ and the two cups of $D_2$. We call such loops
\emph{circular components}, and we do not take account of them
when we define composition of \SJ\ diagrams. When, however, we
define the composition $(D_2,n_2)\cirk(D_1,n_1)=(D_2\cirk
D_1,n_1\pl n_2\pl k)$ of the \SK\ diagrams ${(D_1,n_1)}$ and
${(D_2,n_2)}$, then $k$ is the number of such circular components
that arose as in the picture above. This number  must be finite,
because $D_1$ and $D_2$ have a finite number of cups and caps.

Composition of \SJ\ diagrams, as well as composition of \SK\
diagrams, is associative, and it is clear that the unit \SJ\ and
\SK\ diagrams are unit elements for these compositions; namely we
have

\begin{tabbing}
\hspace{6em}\= (1$\cal J$)\quad\= $I\cirk D=D\cirk I=D$,
\\[1ex]
\> (1$\cal K$)\> $(I,0)\cirk(D,n)=(D,n)\cirk(I,0)=(D,n)$.
\end{tabbing}

To prove strictly the associativity of composition, we would need
a more formal definition of composition, which with the definition
of \SJ\ and \SK\ diagrams we have given here would be rather
involved (cf.\ \cite{DP03c} and \cite{DP03d}). The proof of
associativity would be even more so. A simpler proof of
associativity could be obtained with a topological definition of
diagrams and their composition, on the lines of \cite{DP03a} and
\cite{DP03b} (Section~5), but this other definition would be
rather cumbersome too, and since we have dealt with these matters
in detail in \cite{DP03a} and \cite{DP03b}, we will not go again
into them here.

For ${k\in \N^+}$ let the \emph{cup} \SJ\ diagram $V_k$ be the
\SJ\ diagram that is drawn as follows:

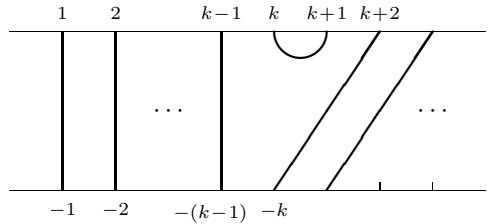
\begin{center}
\begin{picture}(180,100)
{\linethickness{0.05pt} \put(0,20){\line(1,0){180}}
\put(0,80){\line(1,0){180}} \put(140,20){\line(0,1){3}}
\put(160,20){\line(0,1){3}}}

\put(20,15){\makebox(0,0)[t]{\scriptsize $-1$}}
\put(40,15){\makebox(0,0)[t]{\scriptsize $-2$}}
\put(76,15){\makebox(0,0)[t]{\scriptsize $-(k\mn 1)$}}
\put(100,15){\makebox(0,0)[t]{\scriptsize $-k$}}

\put(20,85){\makebox(0,0)[b]{\scriptsize $1$}}
\put(40,85){\makebox(0,0)[b]{\scriptsize $2$}}
\put(80,85){\makebox(0,0)[b]{\scriptsize $k\mn 1$}}
\put(100,85){\makebox(0,0)[b]{\scriptsize $k$}}
\put(120,85){\makebox(0,0)[b]{\scriptsize $k\pl 1$}}
\put(140,85){\makebox(0,0)[b]{\scriptsize $k\pl 2$}}

\thicklines \put(20,20){\line(0,1){60}}
\put(40,20){\line(0,1){60}} \put(80,20){\line(0,1){60}}

\put(100,20){\line(2,3){40}} \put(120,20){\line(2,3){40}}

\put(60,50){\makebox(0,0){$\cdots$}}
\put(160,50){\makebox(0,0){$\cdots$}}

\put(110,80){\oval(20,20)[b]}
\end{picture}
\end{center}

\noindent The \emph{cap} \SJ\ diagram $\Lambda_k$ is the analogous
\SJ\ diagram that is drawn as follows:

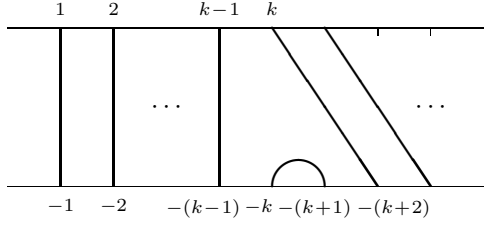
\begin{center}
\begin{picture}(180,100)
{\linethickness{0.05pt} \put(0,20){\line(1,0){180}}
\put(0,80){\line(1,0){180}} \put(140,80){\line(0,-1){3}}
\put(160,80){\line(0,-1){3}}}

\put(20,15){\makebox(0,0)[t]{\scriptsize $-1$}}
\put(40,15){\makebox(0,0)[t]{\scriptsize $-2$}}
\put(74,15){\makebox(0,0)[t]{\scriptsize $-(k\mn 1)$}}
\put(95,15){\makebox(0,0)[t]{\scriptsize $-k$}}
\put(116,15){\makebox(0,0)[t]{\scriptsize $-(k\pl 1)$}}
\put(146,15){\makebox(0,0)[t]{\scriptsize $-(k\pl 2)$}}

\put(20,85){\makebox(0,0)[b]{\scriptsize $1$}}
\put(40,85){\makebox(0,0)[b]{\scriptsize $2$}}
\put(80,85){\makebox(0,0)[b]{\scriptsize $k\mn 1$}}
\put(100,85){\makebox(0,0)[b]{\scriptsize $k$}}

\thicklines \put(20,20){\line(0,1){60}}
\put(40,20){\line(0,1){60}} \put(80,20){\line(0,1){60}}

\put(100,80){\line(2,-3){40}} \put(120,80){\line(2,-3){40}}

\put(60,50){\makebox(0,0){$\cdots$}}
\put(160,50){\makebox(0,0){$\cdots$}}

\put(110,20){\oval(20,20)[t]}
\end{picture}
\end{center}

\noindent The \emph{crossing} \SJ\ diagram $X_k$ is the \SJ\
diagram that is drawn as follows:

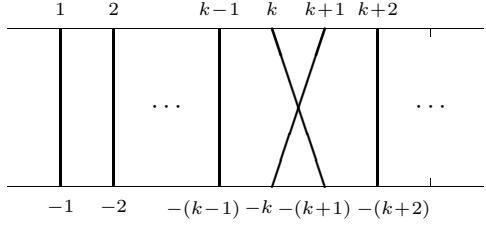
\begin{center}
\begin{picture}(180,100)
{\linethickness{0.05pt} \put(0,20){\line(1,0){180}}
\put(0,80){\line(1,0){180}} \put(160,20){\line(0,1){3}}
\put(160,80){\line(0,-1){3}}}

\put(20,15){\makebox(0,0)[t]{\scriptsize $-1$}}
\put(40,15){\makebox(0,0)[t]{\scriptsize $-2$}}
\put(74,15){\makebox(0,0)[t]{\scriptsize $-(k\mn 1)$}}
\put(95,15){\makebox(0,0)[t]{\scriptsize $-k$}}
\put(116,15){\makebox(0,0)[t]{\scriptsize $-(k\pl 1)$}}
\put(146,15){\makebox(0,0)[t]{\scriptsize $-(k\pl 2)$}}

\put(20,85){\makebox(0,0)[b]{\scriptsize $1$}}
\put(40,85){\makebox(0,0)[b]{\scriptsize $2$}}
\put(80,85){\makebox(0,0)[b]{\scriptsize $k\mn 1$}}
\put(100,85){\makebox(0,0)[b]{\scriptsize $k$}}
\put(120,85){\makebox(0,0)[b]{\scriptsize $k\pl 1$}}
\put(140,85){\makebox(0,0)[b]{\scriptsize $k\pl 2$}}

\thicklines \put(20,20){\line(0,1){60}}
\put(40,20){\line(0,1){60}} \put(80,20){\line(0,1){60}}
\put(140,20){\line(0,1){60}} \put(100,20){\line(1,3){20}}
\put(120,20){\line(-1,3){20}}

\put(60,50){\makebox(0,0){$\cdots$}}
\put(160,50){\makebox(0,0){$\cdots$}}

\end{picture}
\end{center}

\noindent The cup, cap and crossing \SK\ diagrams are ${(V_k,0)}$,
${(\Lambda_k,0)}$ and ${(X_k,0)}$ respectively.

Let \De\ be the set of \SJ\ diagrams. We define inductively as
follows a map $\iota$ from the terms of \SKo\ into \De:

\begin{tabbing}
\hspace{12em}\= $\iota(\dk{k})$ \= = \= $V_k$,
\\[.5ex]
\>$\iota(\gk{k})$ \> = \> $\Lambda_k$,
\\[.5ex]
\> $\iota(\sigma_k)$ \> = \> $X_k$,
\\[.5ex]
\> $\iota(\mj)$ \> = \> $I$,
\\[.5ex]
\> $\iota(tu)$ \> = \> $\iota(t)\cirk\iota(u)$.
\end{tabbing}

\noindent By this definition, we have that

\[
\iota(c)=\iota(\dk{k}\gk{k})=\iota(\dk{k})\cirk\iota(\gk{k})=V_k\cirk\Lambda_k=I.
\]

The map $\kappa$ from the terms of \SKo\ into ${\De\times N}$ is
defined by

\begin{tabbing}
\hspace{5em}\= $\kappa(tu)$ \= = \= \kill

\> $\kappa(t)$ \> = \> $(\iota(t),0)$,\quad where $t$ is $\dk{k}$,
$\gk{k}$, $\sigma_k$ or \mj,
\\*[1ex]
\> $\kappa(tu)$ \> = \> $\kappa(t)\cirk\kappa(u)$.
\end{tabbing}

\noindent By this definition we have that

\[
\kappa(c)=\kappa(\dk{k}\gk{k})=\kappa(\dk{k})\cirk\kappa(\gk{k})=(V_k,0)\cirk(\Lambda_k,0)=(I,1).
\]

Then we can prove the following lemma.

\prop{Soundness Lemma}{\\[1ex] \indent$(\cal K)$\quad If $t=u$ in \SKo,
then $\kappa(t)=\kappa(u)$. \\[1ex]\indent
$(\cal J)$\quad If $t=u$ in \SJo, then $\iota(t)=\iota(u)$.}

\dkz We have mentioned above that we have (1$\cal J$) and (1$\cal
K$), and that $\cal J$ and $\cal K$ composition is associative,
which takes care of the equations (1) and (2). It is
straightforward to verify the remaining equations of \SKo\ and
\SJo. \qed

\vspace{2ex}

This lemma says that $\kappa$ and $\iota$ give rise to
homomorphisms from the monoids \SKo\ and \SJo\ to monoids
${(\De\times N,\cirk,(I,0))}$ and ${(\De,\cirk,I)}$ respectively.
Our goal in the following sections is to show that these
homomorphisms are isomorphisms.

\section{Normal forms in \SKo\ and \SJo}

Let ${i,j\in \N^+}$, and let ${j\leq i}$. We define as follows the
term $\dk{j,i}$ of \SKo, which we call a \emph{block-cup}:
$\dk{j,j}$ is $\dk{j}$, and if ${j<i}$, then $\dk{j,i}$ is
${\dk{j,i\mn 1}\sigma_i}$.

We define analogously the \emph{block-cap} $\gk{i,j}$: we have
that $\gk{j,j}$ is $\gk{j}$, and if ${j<i}$, then $\gk{i,j}$ is
${\sigma_i\gk{i\mn 1,j}}$.

For $i$ and $j$ as above, let the \emph{block-crossing} be the
term $[i,j]$ of \SKo\ defined as follows: $[j,j]$ is $\sigma_j$,
and if ${j<i}$, then ${[i,j]}$ is ${\sigma_i[i\mn 1,j]}$.

It is easy to see that for ${j<i}$ the \SJ\ diagram
${\iota(\dk{j,i})}$, which we call the \emph{block-cup} \SJ\
\emph{diagram}, is drawn as follows:

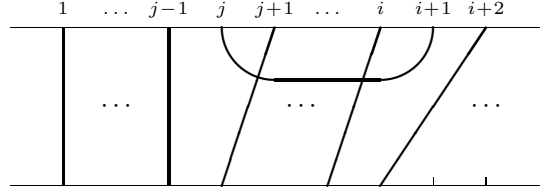
\begin{center}
\begin{picture}(200,100)
{\linethickness{0.05pt} \put(0,20){\line(1,0){200}}
\put(0,80){\line(1,0){200}} \put(160,20){\line(0,1){3}}
\put(180,20){\line(0,1){3}}}

\put(20,85){\makebox(0,0)[b]{\scriptsize $1$}}
\put(40,85){\makebox(0,0)[b]{\scriptsize $\ldots$}}
\put(60,85){\makebox(0,0)[b]{\scriptsize $j\mn 1$}}
\put(80,85){\makebox(0,0)[b]{\scriptsize $j$}}
\put(100,85){\makebox(0,0)[b]{\scriptsize $j\pl 1$}}
\put(120,85){\makebox(0,0)[b]{\scriptsize $\ldots$}}
\put(140,85){\makebox(0,0)[b]{\scriptsize $i$}}
\put(160,85){\makebox(0,0)[b]{\scriptsize $i\pl 1$}}
\put(180,85){\makebox(0,0)[b]{\scriptsize $i\pl 2$}}

\thicklines \put(20,20){\line(0,1){60}}
\put(60,20){\line(0,1){60}} \put(80,20){\line(1,3){20}}
\put(120,20){\line(1,3){20}} \put(140,20){\line(2,3){40}}

\put(40,50){\makebox(0,0){$\cdots$}}
\put(110,50){\makebox(0,0){$\cdots$}}
\put(180,50){\makebox(0,0){$\cdots$}}

\put(120,80){\oval(80,40)[b]}
\end{picture}
\end{center}

\noindent Analogously, the \emph{block-cap} \SJ\ \emph{diagram}
${\iota(\gk{i,j})}$ is drawn as follows:

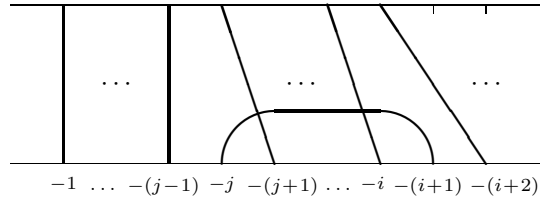
\begin{center}
\begin{picture}(200,100)
{\linethickness{0.05pt} \put(0,20){\line(1,0){200}}
\put(0,80){\line(1,0){200}} \put(160,80){\line(0,-1){3}}
\put(180,80){\line(0,-1){3}}}

\put(20,15){\makebox(0,0)[t]{\scriptsize $-1$}}
\put(35,11){\makebox(0,0)[t]{\scriptsize $\ldots$}}
\put(58,15){\makebox(0,0)[t]{\scriptsize $-(j\mn 1)$}}
\put(80,15){\makebox(0,0)[t]{\scriptsize $-j$}}
\put(103,15){\makebox(0,0)[t]{\scriptsize $-(j\pl 1)$}}
\put(124,11){\makebox(0,0)[t]{\scriptsize $\ldots$}}
\put(137,15){\makebox(0,0)[t]{\scriptsize $-i$}}
\put(158,15){\makebox(0,0)[t]{\scriptsize $-(i\pl 1)$}}
\put(187,15){\makebox(0,0)[t]{\scriptsize $-(i\pl 2)$}}

\thicklines \put(20,20){\line(0,1){60}}
\put(60,20){\line(0,1){60}} \put(100,20){\line(-1,3){20}}
\put(140,20){\line(-1,3){20}} \put(180,20){\line(-2,3){40}}

\put(40,50){\makebox(0,0){$\cdots$}}
\put(110,50){\makebox(0,0){$\cdots$}}
\put(180,50){\makebox(0,0){$\cdots$}}

\put(120,20){\oval(80,40)[t]}
\end{picture}
\end{center}

\noindent And the \emph{block-crossing} \SJ\ \emph{diagram}
${\iota([i,j])}$ is drawn as follows:

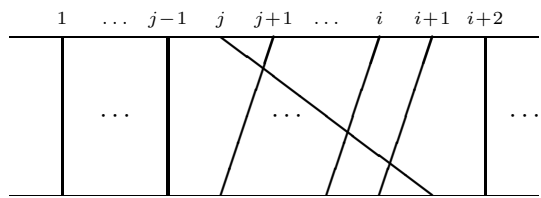
\begin{center}
\begin{picture}(200,100)
{\linethickness{0.05pt} \put(0,20){\line(1,0){200}}
\put(0,80){\line(1,0){200}}}

\put(20,85){\makebox(0,0)[b]{\scriptsize $1$}}
\put(40,85){\makebox(0,0)[b]{\scriptsize $\ldots$}}
\put(60,85){\makebox(0,0)[b]{\scriptsize $j\mn 1$}}
\put(80,85){\makebox(0,0)[b]{\scriptsize $j$}}
\put(100,85){\makebox(0,0)[b]{\scriptsize $j\pl 1$}}
\put(120,85){\makebox(0,0)[b]{\scriptsize $\ldots$}}
\put(140,85){\makebox(0,0)[b]{\scriptsize $i$}}
\put(160,85){\makebox(0,0)[b]{\scriptsize $i\pl 1$}}
\put(180,85){\makebox(0,0)[b]{\scriptsize $i\pl 2$}}

\thicklines \put(20,20){\line(0,1){60}}
\put(60,20){\line(0,1){60}} \put(80,20){\line(1,3){20}}
\put(120,20){\line(1,3){20}} \put(140,20){\line(1,3){20}}
\put(160,20){\line(-4,3){80}} \put(180,20){\line(0,1){60}}

\put(40,50){\makebox(0,0){$\cdots$}}
\put(105,50){\makebox(0,0){$\cdots$}}
\put(195,50){\makebox(0,0){$\cdots$}}

\end{picture}
\end{center}

Consider terms of \SKo\ of the form

\[
c^l\gk{i_p,j_p}\ldots\gk{i_1,j_1}\mj[k_1,l_1]\ldots[k_q,l_q]\dk{m_1,n_1}\ldots\dk{m_r,n_r}
\]

\noindent where ${l,p,q,r\geq 0}$, ${j_1<\ldots<j_p}$,
${k_1<\ldots<k_q}$ and ${m_1<\ldots<m_r}$. The notation $c^l$ has
been introduced in Section~2. If $p=0$, then the sequence
${\gk{i_p,j_p}\ldots\gk{i_1,j_1}}$ is empty, and analogously if
$q=0$, or $r=0$. Terms of this form will be said to be in
\emph{normal form}.

We can then prove the following lemma for \SKo.

\prop{Normal Form Lemma}{Every term is equal to a term in normal
form.}

\dkz We will give a reduction procedure that transforms every term
$t$ of \SKo\ into a term $t'$ in normal form such that ${t=t'}$ in
\SKo. This procedure will have three phases, but it might
terminate after any of the phases (or before any phase, if the
term we start from is in normal form).

Roughly speaking, in the first phase of our procedure we push caps
to the left, where they are transformed into block-caps, which we
arrange in the decreasing order of the normal form. Circles we
might have at the beginning and circles that arise during the
first phase are all pushed to the left in this phase. Similarly,
unit terms \mj\ that we might have at the beginning and that arise
during the first phase are all deleted, except if the whole term
is \mj. If the whole term is not \mj, we put one term \mj\ at an
appropriate place. This is the first phase of the procedure.

In the second phase we push cups to the right, where they are
transformed into block-cups, which we arrange in the increasing
order of the normal form. (We could as well have decided to push
cups to the right in the first phase, and then caps to the left in
the second phase.)

After the first and second phase are over, the crossings that may
remain in between the block-caps on the left and the block-cups on
the right are transformed into block-crossings, which we arrange
in the increasing order of the normal form. This is the third, and
last, phase of the procedure.

At the beginning of the first phase of the reduction procedure we
replace every cap $\gk{i}$ by the block-cap $\gk{i,i}$. Then,
according to the following equations of \SKo, which are derived in
a straightforward manner, we replace every subterm of the form on
the left-hand side by the term on the right-hand side of these
equations:

\begin{tabbing}
\hspace{1em}\emph{cup-block-cap equations}:
\\*[1.5ex]
\hspace{2em}\= (\emph{cup-cap}~i)\quad \= $\dk{k}\gk{i,j}=\gk{i\mn
2,j\mn 2}\dk{k}$\quad if $k\pl 2\leq j$,
\\[1ex]
\> (\emph{cup-cap}~ii) \> for $k\leq j\leq k\pl 1$,
\\*[.5ex]
\hspace{3em}\= (\emph{cup-cap}~ii.1)\quad \=
$\dk{k}\gk{i,j}=\sigma_{i-2}\ldots\sigma_k$\quad \= if $k\pl 2\leq
i$,
\\*[.5ex]
\> (\emph{cup-cap}~ii.2)\> $\dk{k}\gk{i,j}=\mj$\> if $i=k\pl 1$,
\\*[.5ex]
\> (\emph{cup-cap}~ii.3)\> $\dk{k}\gk{i,j}=c$\> if $i=k$,
\\[1ex]
\hspace{2em}\= (\emph{cup-cap}~iii)\quad \= for $j\leq k\mn 1$,
\\*[.5ex]
\hspace{3em}\= (\emph{cup-cap}~iii.1)\quad \=
$\dk{k}\gk{i,j}=\gk{i\mn 2,j}\dk{k\mn 1}$\quad if $k\leq i\mn 1$,
\\[.5ex]
\> (\emph{cup-cap}~iii.2)\> for $ i\leq k\leq i\pl 1$,
\\*[.5ex]
\hspace{4em}\= (\emph{cup-cap}~iii.2.1)\quad \=
$\dk{k}\gk{i,j}=\sigma_j\ldots\sigma_{k-2}$\quad \= if $j<k\mn 1$,
\\*[.5ex]
\> (\emph{cup-cap}~iii.2.2)\> $\dk{k}\gk{i,j}=\mj$\> if $j=k\mn
1$,
\\[.5ex]
\hspace{3em}\= (\emph{cup-cap}~iii.1)\quad \=\kill

\> (\emph{cup-cap}~iii.3)\> $\dk{k}\gk{i,j}=\gk{i,j}\dk{k\mn 2}$\>
if $i\pl 2\leq k$;
\\[1.5ex]
\hspace{1em}$\mj$ \emph{and} $c$ \emph{equations}:
\\*[1.5ex]
\hspace{2em}\= (1)\quad \= $\mj t=t$,\quad $t\mj=t$,
\\*[.5ex]
\> ($c$)\> $tc=ct$;
\\[1.5ex]
\hspace{1em}$\sigma$-\emph{block-cap equations}:
\\*[1.5ex]
\hspace{2em}\= ($\sigma$-\emph{cap}~i)\quad \=
$\sigma_k\gk{i,j}=\gk{i,j}\sigma_k$\quad \= if $k\pl 2\leq j$,
\\[1ex]
\> ($\sigma$-\emph{cap}~ii)\> $\sigma_k\gk{i,j}=\gk{i,j\mn 1}$\>
if $j=k\pl 1$,
\\[1ex]
\> ($\sigma$-\emph{cap}~iii)\> for $j=k$,
\\*[.5ex]
\hspace{3em}\= ($\sigma$-\emph{cap}~iii.1)\quad \=
$\sigma_k\gk{i,j}=\gk{i,j\pl 1}$\quad \= if $j<i$,
\\*[.5ex]
\> ($\sigma$-\emph{cap}~iii.2)\> $\sigma_k\gk{i,j}=\gk{i,j}$\> if
$i=j$,
\\[1ex]
\hspace{2em}\= ($\sigma$-\emph{cap}~i)\quad \= \kill

\> ($\sigma$-\emph{cap}~iv)\> for $j\leq k\mn 1$,
\\*[.5ex]
\hspace{3em}\= ($\sigma$-\emph{cap}~iv.1)\quad \=
$\sigma_k\gk{i,j}=\gk{i,j}\sigma_{k-1}$\quad \= if $k\pl 1\leq i$,
\\*[.5ex]
\> ($\sigma$-\emph{cap}~iv.2)\> $\sigma_k\gk{i,j}=\gk{i\mn 1,j}$\>
if $i=k$,
\\*[.5ex]
\> ($\sigma$-\emph{cap}~iv.3)\> $\sigma_k\gk{i,j}=\gk{i\pl 1,j}$\>
if $i=k\mn 1$,
\\*[.5ex]
\> ($\sigma$-\emph{cap}~iv.4)\>
$\sigma_k\gk{i,j}=\gk{i,j}\sigma_{k-2}$\> if $i\leq k\mn 2$ (this
case is impossible
\\
\` if $j=k\mn 1$);
\\[1.5ex]
\hspace{1em}\emph{double block-cap equations}:
\\*[1.5ex]
\hspace{3em} for $j\leq l$
\\*[1ex]
\hspace{2em}\= (\emph{cap}~i)\quad \= $\gk{i,j}\gk{k,l}=\gk{k\pl
2,l\pl 2}\gk{i,j}$\quad\quad\quad\= if $i\leq l$,
\\*[1ex]
\> (\emph{cap}~ii)\> $\gk{i,j}\gk{k,l}=\gk{k\pl 2,l\pl 1}\gk{i\mn
1,j}$\> if $l\pl 1\leq i\leq k\pl 1$,
\\*[1ex]
\> (\emph{cap}~iii)\> $\gk{i,j}\gk{k,l}=\gk{k\pl 1,l\pl 1}\gk{i\mn
2,j}$\> if $k\pl 2\leq i$.
\end{tabbing}

\noindent Note that with the equation ($c$) we can push all
circles to the extreme left. Note next that the cup-block-cap
equations cover all cases for subterms of the form
${\dk{k}\gk{i,j}}$, and the $\sigma$-block-cap equations cover all
cases for subterms of the form ${\sigma_k\gk{i,j}}$. The double
block-cap equations cover all cases for subterms of the form
${\gk{i,j}\gk{k,l}}$ with ${j\leq l}$, and they replace such
subterms by ${\gk{k',l'}\gk{i',j'}}$ with ${l'>j'}$. So, when we
can no longer make any replacement by one of the equations above,
we are left with the term $\mj$ or with a term of the form

\[
c^l\gk{i_p,j_p}\ldots\gk{i_1,j_1}u_1\ldots u_n
\]

\noindent where ${l,p,n\geq 0}$, ${l\pl p\pl n\geq 1}$,
${j_1<\ldots < j_p}$, and the term ${u_1\ldots u_n}$, which is the
empty sequence if ${n=0}$, is made exclusively of cup terms and
crossing terms. In case we are not left with $\mj$, we replace the
term above with

\[
c^l\gk{i_p,j_p}\ldots\gk{i_1,j_1}\mj u_1\ldots u_n.
\]

\noindent With that the first phase of the reduction procedure is
over.

If ${n=0}$, then the whole procedure is over, and we are left with
a term in normal form. If ${n\geq 1}$, then we enter into the
second phase of the procedure, and reduce ${u_1\ldots u_n}$. We
first replace every cup $\dk{i}$ in this term by the block-cup
${\dk{i,i}}$, and then we make replacements of subterms of the
form ${\dk{j,i}\sigma_k}$ and ${\dk{l,k}\dk{j,i}}$ with ${j\leq
l}$ by using equations exactly dual to the $\sigma$-block-cap
equations and the double block-cap equations. When we can no
longer make any of these replacements, we have instead of
${u_1\ldots u_n}$ a term of the form

\[
v_1\ldots v_m\dk{m_1,n_1}\ldots\dk{m_r,n_r}
\]

\noindent where ${m,r\geq 0}$, ${m\pl r\geq 1}$,
${m_1<\ldots<m_r}$, and the term ${v_1\ldots v_m}$, which is the
empty sequence if ${m=0}$, is made exclusively of crossing terms.
With that the second phase of the reduction procedure is over.

If ${m=0}$, then the whole reduction procedure is over, and we are
left with a term in normal form. If ${m\geq 1}$, then we enter
into the third phase of the procedure, and reduce ${v_1\ldots
v_m}$. We first replace every crossing term $\sigma_i$ in this
term by the block-crossing ${[i,i]}$. Then, according to the
following equations of \SKo, which are derived in a
straightforward manner from the monoid equations and the
$\sigma$~equations of Section~2, we replace every subterm of the
form on the left-hand side by the term on the right-hand side of
these equations:

\begin{tabbing}
\hspace{1em}\emph{block}-$\sigma$ \emph{equations}:
\\*[1.5ex]
\hspace{3em} for $k\leq i$
\\*[.5ex]
\hspace{2em}\= ($\sigma$i)\quad \= $[i,j][k,l]=[k,l][i,j]$\quad \=
if $k\pl 2\leq j$,
\\[1ex]
\> ($\sigma$ii)\> $[i,j][k,l]=[i,l]$\> if $j=k\pl 1$,
\\[1ex]
\> ($\sigma$iii)\> for $j=k$,
\\*[.5ex]
\hspace{3em}\= ($\sigma$iii.1)\quad \= $[i,j][k,l]=[k\mn
1,l][i,j\pl 1]$\quad \= if $j<i$, $l<k$,
\\*[.5ex]
\> ($\sigma$iii.2)\> $[i,j][k,l]=[i,j\pl 1]$\> if $j<i$, $l=k$,
\\*[.5ex]
\> ($\sigma$iii.3)\> $[i,j][k,l]=[k\mn 1,l]$\> if $j=i$, $l<k$,
\\*[.5ex]
\> ($\sigma$iii.4)\> $[i,j][k,l]=\mj$\> if $j=i$, $l=k$,
\\[1ex]
\hspace{2em}\= ($\sigma$i)\quad \= \hspace{2em}\= \kill

\> ($\sigma$iv)\> for $j\leq k\mn 1$,
\\*[.5ex]
\hspace{3em}\= ($\sigma$iv.1)\quad \= $[i,j][k,l]=[k\mn
1,l][i,j\pl 1]$\quad \= if $l\leq j$,
\\*[.5ex]
\> ($\sigma$iv.2)\> $[i,j][k,l]=[k\mn 1,l\mn 1][i,j]$\> if $j<l$;
\\[1.5ex]
\hspace{1em}$\mj$ \emph{equations}:
\\*[1.5ex]
\hspace{2em}\= (1)\quad \= $\mj t=t$,\quad $t\mj=t$.
\end{tabbing}

Note that the block-$\sigma$ equations cover all cases for
subterms of the form ${[i,j][k,l]}$ with ${k\leq i}$, and they
replace such subterms by ${[k',l'][i',j']}$ with ${k'<i'}$, or by
single block-crossings, or by $\mj$. The terms $\mj$ are
eliminated by the $\mj$ equations, except if the whole term is
$\mj$. When we can no longer make any of these replacements, we
have instead of ${v_1\ldots v_m}$ either $\mj$ or a term of the
form

\[
[k_1,l_1]\ldots[k_q,l_q]
\]

\noindent where ${q\geq 1}$ and ${k_1<\ldots<k_q}$. With that the
third phase of the reduction procedure is over. The original term
has now been reduced to normal form, except for one more possible
application of the $\mj$ equations if ${v_1\ldots v_m}$ has been
replaced by $\mj$ and ${l\pl p\pl r\geq 1}$. \qed

\vspace{2ex}

Note that the cup-cap equations of Section~2 are all instances of
the block equations mentioned in the course of this proof. The
equation (\emph{cap}) is an instance of (\emph{cap}~i) for ${i=j}$
and ${k=l}$, and analogously for (\emph{cup}) and the dual
equation with block-cups, which we mentioned in the proof, but did
not write down. The equation (\emph{cup-cap}~1) is an instance of
(\emph{cup-cap}~iii.3) for ${i=j}$, the equation
(\emph{cap-cup}~1) is an instance of (\emph{cup-cap}~i) for
${i=j}$, and the equations (\emph{cup-cap}) are instances of
(\emph{cup-cap}~ii.2) for ${i=j}$ and of (\emph{cup-cap}~iii.2.2)
for ${i=j}$.

The same holds for the $\sigma$-cup and $\sigma$-cap equations of
Section~2. The equation ($\sigma$-\emph{cap}~1) is an instance of
($\sigma$-\emph{cap}~i) for ${i=j}$, the equation
($\sigma$-\emph{cap}~2) is an instance of
($\sigma$-\emph{cap}~iv.4) for ${i=j}$, the equation
($\sigma$-\emph{cap}~3) is ($\sigma$-\emph{cap}~iii.2), and
($\sigma$-\emph{cap}~4) is an instance of ($\sigma$-\emph{cap}~ii)
for ${i=j}$. We have an analogous situation for the $\sigma$-cup
equations with respect to the equations dual to the
$\sigma$-block-cap equations with block-cups, which we mentioned
in the proof, but did not write down.

Finally, the $\sigma$ equations of Section~2 are instances of the
block-$\sigma$ equations of the proof above. The equation
($\sigma$) is an instance of ($\sigma$i) for ${i=j}$ and ${k=l}$,
the equation ($\sigma$2) is ($\sigma$iii.4), and ($\sigma$3) is an
instance of ($\sigma$iv.2) for ${i=j\pl 1=k=l}$.

So if for presenting \SKo\ we had taken as generators block-cups,
block-caps and block-crossings, instead of cups, caps and
crossings, the block equations of the reduction procedure of our
proof of the Normal Form Lemma would suffice for this
presentation, though some equations are redundant. (In this block
presentation, the equations ($\sigma$-\emph{cap}~iv.3) and
($\sigma$ii) are, however essential, in order to justify the
definition of block-caps and block-crossings in terms of caps and
crossings.)

Note that if we restrict ourselves to the third phase of the
proof, we have provided a normal form for terms of symmetric
groups, and an alternative presentation, with block-crossings as
generators, for these groups. This normal form for symmetric
groups, which may be  derived from \cite{B11} (Note~C, pp.\
464-469), is analogous up to a point to the Jones normal form for
terms of the monoids of Temperley-Lieb algebras (see \cite{DP03a}
and \cite{DP03b}, Section~9). As we shall see later (Section~6),
this normal form serves well to show the completeness of the
presentation of symmetric groups with the $\sigma$ equations of
Section~2, or with the block-$\sigma$ equations of our proof of
the Normal Form Lemma. (This completeness proof is different that
in \cite{B11}, Note~C; cf.\ \cite{DP04}, Section 5.2.)

If in the definition of normal form for terms of \SKo\ we require
that ${l=0}$, i.e.\ that $c^l$ is the empty sequence (there are no
circles at the beginning), we obtain the definition of normal form
appropriate for \SJo, which we call $\cal J$-\emph{normal form}.
Then our proof of the Normal Form Lemma suffices to prove this
lemma for \SJo\ with respect to the $\cal J$-normal form. We just
may have to apply at the end of our procedure the equation ($c$1)
of Section~2, together with the equation (1), to get rid of
circles.

\section{Generating \SK\ and \SJ\ diagrams}

A cup ${\{j,i\}}$ of an \SJ\ diagram $D$ is \emph{maximal} when
there is no cup ${\{j',i'\}}$ of $D$ such that ${j<j'}$. (Remember
that according to our convention of Section~3 we have here ${j<i}$
and ${j'<i'}$.) Analogously, a cap ${\{-j,-i\}}$ is \emph{maximal}
when there is no cap ${\{-j',-i'\}}$ of $D$ such that ${j<j'}$. If
there are cups in $D$, then there must be a maximal one because
there are only finitely many cups in $D$, and analogously for
caps.

Let $P$ be an \SJ\ diagram without cups and caps. Transversal
threads ${\{l,-k\}}$ of $P$ such that ${l<k}$ will be called
\emph{falling}. A falling thread of $P$ is \emph{maximal} when
there is no falling thread ${\{l',k'\}}$ of $P$ such that
${k<k'}$. If there are falling threads in $P$, then there must be
a maximal one, because $P$ is of type ${(n,n)}$ for some ${n\in
\N^+}$.

Then we can prove the following lemma.

\prop{Generating Lemma}{Every \SJ\ diagram is equal to a diagram
obtained from cup \SJ\ diagrams, cap \SJ\ diagrams, crossing \SJ\
diagrams and the unit \SJ\ diagram with the help of composition.}

\dkz Let $D$ be an \SJ\ diagram. If there are cups in $D$, then
there must be a maximal cup ${\{m,n\pl 1\}}$, with ${m\leq n}$.
Then consider the partition of ${\Z-\{0,m,n\pl 1\}}$ obtained from
$D$ just by omitting the cup ${\{m,n\pl 1\}}$. This partition
gives rise to an \SJ\ diagram $D'$ by replacing every ${x\in \Z}$
such that ${m\pl 1\leq x\leq n}$ by ${x\mn 1}$, and every ${y\in
\Z}$ such that ${n\pl 2\leq y}$ by ${y\mn 2}$. If $D$ is of type
${(s_1,s_2)}$, then $D'$ is of type ${(s_1\mn 2,s_2)}$. Then we
can check that ${D=D'\cirk \iota(\dk{m,n})}$, and $D'$ has one cup
less than $D$. For example, we have

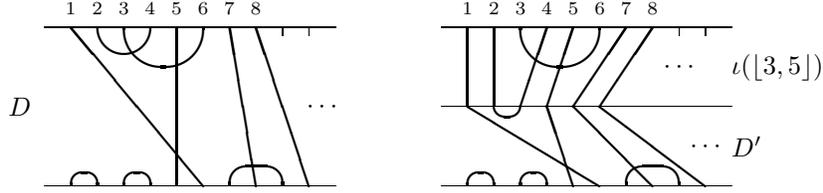
\begin{center}
\begin{picture}(280,100)
{\linethickness{0.05pt} \put(0,20){\line(1,0){110}}
\put(0,80){\line(1,0){110}} \put(90,80){\line(0,-1){3}}
\put(100,80){\line(0,-1){3}} \put(150,20){\line(1,0){110}}
\put(150,50){\line(1,0){110}} \put(150,80){\line(1,0){110}}
\put(240,80){\line(0,-1){3}} \put(250,80){\line(0,-1){3}}}

\put(10,85){\makebox(0,0)[b]{\scriptsize $1$}}
\put(20,85){\makebox(0,0)[b]{\scriptsize $2$}}
\put(30,85){\makebox(0,0)[b]{\scriptsize $3$}}
\put(40,85){\makebox(0,0)[b]{\scriptsize $4$}}
\put(50,85){\makebox(0,0)[b]{\scriptsize $5$}}
\put(60,85){\makebox(0,0)[b]{\scriptsize $6$}}
\put(70,85){\makebox(0,0)[b]{\scriptsize $7$}}
\put(80,85){\makebox(0,0)[b]{\scriptsize $8$}}

\put(160,85){\makebox(0,0)[b]{\scriptsize $1$}}
\put(170,85){\makebox(0,0)[b]{\scriptsize $2$}}
\put(180,85){\makebox(0,0)[b]{\scriptsize $3$}}
\put(190,85){\makebox(0,0)[b]{\scriptsize $4$}}
\put(200,85){\makebox(0,0)[b]{\scriptsize $5$}}
\put(210,85){\makebox(0,0)[b]{\scriptsize $6$}}
\put(220,85){\makebox(0,0)[b]{\scriptsize $7$}}
\put(230,85){\makebox(0,0)[b]{\scriptsize $8$}}

\put(260,35){\makebox(0,0)[l]{$D'$}}
\put(260,65){\makebox(0,0)[l]{$\iota(\dk{3,5})$}}

\thicklines \put(50,20){\line(0,1){60}}
\put(60,20){\line(-5,6){50}} \put(80,20){\line(-1,6){10}}
\put(100,20){\line(-1,3){20}}

\put(160,50){\line(0,1){30}} \put(170,50){\line(0,1){30}}
\put(180,50){\line(1,3){10}} \put(190,50){\line(1,3){10}}
\put(200,50){\line(2,3){20}} \put(210,50){\line(2,3){20}}

\put(200,20){\line(-1,3){10}} \put(210,20){\line(-5,3){50}}
\put(230,20){\line(-1,1){30}} \put(250,20){\line(-4,3){40}}

\put(105,50){\makebox(0,0){$\cdots$}}
\put(250,35){\makebox(0,0){$\cdots$}}
\put(240,65){\makebox(0,0){$\cdots$}}

\put(15,20){\oval(10,10)[t]} \put(35,20){\oval(10,10)[t]}
\put(80,20){\oval(20,15)[t]}

\put(30,80){\oval(20,20)[b]} \put(45,80){\oval(30,30)[b]}

\put(165,20){\oval(10,10)[t]} \put(185,20){\oval(10,10)[t]}
\put(230,20){\oval(20,15)[t]}

\put(175,50){\oval(10,8)[b]} \put(195,80){\oval(30,30)[b]}

\put(-5,50){\makebox(0,0)[r]{$D$}}

\end{picture}
\end{center}

We proceed analogously, pulling out maximal caps on the left, if
$D$ has caps.

By repeated applications of this procedure we must end up with an
\SJ\ diagram equal to $D$ of the form

\[
\iota(\gk{i_p,j_p})\cirk\ldots\cirk\iota(\gk{i_1,j_1})\cirk
P\cirk\iota(\dk{m_1,n_1})\cirk\ldots\cirk\iota(\dk{m_r,n_r})
\]

\noindent where ${p,r\geq 0}$, ${p\pl r\geq 0}$,
${j_1<\ldots<j_p}$, ${m_1<\ldots<m_r}$, and there are no cups and
caps in the \SJ\ diagram $P$. The \SJ\ diagram $P$ corresponds to
a permutation of an initial segment of $\N^+$, and it has a type
${(n,n)}$ for some ${n\in \N^+}$.

If there are no falling threads in $P$, then $P$ is equal to the
unit \SJ\ diagram $I$. If there are falling threads in $P$, then
let ${\{l_q,-(k_q\pl 1)\}}$ with ${l_q\leq k_q}$ be the maximal
falling thread of $P$. Then consider the partition of
${\Z-\{0,l_q,-(k_q\pl 1)\}}$ obtained from $P$ just by omitting
the thread ${\{l_q,-(k_q\pl 1)\}}$. This partition gives rise to
an \SJ\ diagram $P'$ without cups and caps by replacing every
${x\in \Z}$ such that ${l_q+1\leq x\leq k_q+1}$ by ${x\mn 1}$, and
by adding the thread ${\{k_q\pl 1,-(k_q\pl 1)\}}$. If there are
falling threads in $P'$ (i.e., if $P'$ is not the unit \SJ\
diagram $I$), then its maximal falling thread ${\{l,-(k\pl 1)\}}$
must have ${k<k_q}$ (otherwise, ${\{l_q,-(k_q\pl 1)\}}$ would not
be maximal in $P$). Then we can check that
${P=P'\cirk\iota([k_q,l_q])}$. For example, we have

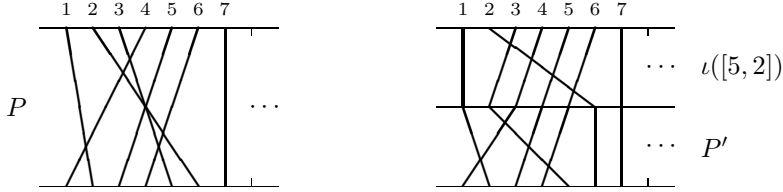
\begin{center}
\begin{picture}(260,100)
{\linethickness{0.05pt} \put(0,20){\line(1,0){90}}
\put(0,80){\line(1,0){90}} \put(150,20){\line(1,0){90}}
\put(150,50){\line(1,0){90}} \put(150,80){\line(1,0){90}}
\put(80,80){\line(0,-1){3}} \put(80,20){\line(0,1){3}}
\put(230,80){\line(0,-1){3}} \put(230,20){\line(0,1){3}}}

\put(10,85){\makebox(0,0)[b]{\scriptsize $1$}}
\put(20,85){\makebox(0,0)[b]{\scriptsize $2$}}
\put(30,85){\makebox(0,0)[b]{\scriptsize $3$}}
\put(40,85){\makebox(0,0)[b]{\scriptsize $4$}}
\put(50,85){\makebox(0,0)[b]{\scriptsize $5$}}
\put(60,85){\makebox(0,0)[b]{\scriptsize $6$}}
\put(70,85){\makebox(0,0)[b]{\scriptsize $7$}}

\put(160,85){\makebox(0,0)[b]{\scriptsize $1$}}
\put(170,85){\makebox(0,0)[b]{\scriptsize $2$}}
\put(180,85){\makebox(0,0)[b]{\scriptsize $3$}}
\put(190,85){\makebox(0,0)[b]{\scriptsize $4$}}
\put(200,85){\makebox(0,0)[b]{\scriptsize $5$}}
\put(210,85){\makebox(0,0)[b]{\scriptsize $6$}}
\put(220,85){\makebox(0,0)[b]{\scriptsize $7$}}

\put(250,35){\makebox(0,0)[l]{$P'$}}
\put(250,65){\makebox(0,0)[l]{$\iota([5,2])$}}

\thicklines \put(10,20){\line(1,2){30}}
\put(20,20){\line(-1,6){10}} \put(30,20){\line(1,3){20}}
\put(40,20){\line(1,3){20}} \put(50,20){\line(-1,3){20}}
\put(60,20){\line(-2,3){40}} \put(70,20){\line(0,1){60}}

\put(160,50){\line(0,1){30}} \put(170,50){\line(1,3){10}}
\put(180,50){\line(1,3){10}} \put(190,50){\line(1,3){10}}
\put(200,50){\line(1,3){10}} \put(210,50){\line(-4,3){40}}
\put(220,50){\line(0,1){30}}

\put(160,20){\line(2,3){20}} \put(170,20){\line(-1,3){10}}
\put(180,20){\line(1,3){10}} \put(190,20){\line(1,3){10}}
\put(200,20){\line(-1,1){30}} \put(210,20){\line(0,1){30}}
\put(220,20){\line(0,1){30}}

\put(85,50){\makebox(0,0){$\cdots$}}
\put(235,35){\makebox(0,0){$\cdots$}}
\put(235,65){\makebox(0,0){$\cdots$}}

\put(-5,50){\makebox(0,0)[r]{$P$}}
\end{picture}
\end{center}

By repeated applications of this procedure we must end up with an
\SJ\ diagram equal to $P$ of the form

\[
I\cirk\iota([k_1,l_1])\cirk\ldots\cirk\iota([k_q,l_q])
\]

\noindent where ${q\geq 0}$ and ${k_1<\ldots<k_q}$.

By definition, block-cup, block-cap and block-crossing \SJ\
diagrams are generated by cup \SJ\ diagrams, cap \SJ\ diagrams and
crossing \SJ\ diagrams with the help of composition. So the lemma
is proved. \qed

\vspace{2ex}

The last part of the proof, which involves block-crossing \SJ\
diagrams, is about a well-known fact. We have, nevertheless,
preferred to go through it in order to show that it is of the same
sort as the previous parts of the proof. This part of the proof is
also a part of our proof of the completeness of the standard
presentation of symmetric groups.

As a matter of fact, we need not have strived in the proof above
to get ${j_1<\ldots<j_p}$, ${n_1<\ldots<n_r}$ and
${k_1<\ldots<k_q}$. However, if we do so, we shall end up with a
composition of \SJ\ diagrams corresponding to a term of \SKo\ in
$\cal J$-normal form.

Strictly speaking, unit \SJ\ diagrams can be omitted in the
formulation of the Generating Lemma since $I=V_k\cirk\Lambda_{k+
1}=V_{k+ 1}\cirk\Lambda_k=X_k\cirk X_k$ for any ${k\in \N^+}$.

The Generating Lemma is provable also when in its formulation
$\cal J$ is replaced everywhere by $\cal K$. For any \SK\ diagram
${(D,n)}$ we have

\[
(D,n)=(I,n)\cirk(D,0).
\]

\noindent The \SK\ diagram ${(I,n)}$ can be generated from cup
\SK\ diagrams, cap \SK\ diagrams and the unit \SK\ diagram with
the help of composition since we have

\[
(I,n\pl 1)=(I,n)\cirk(I,1)=(I,n)\cirk(V_k,0)\cirk(\Lambda_k,0).
\]

\noindent For ${(D,0)}$, we first go through our proof of the
Generating Lemma for \SJ\ diagrams in order to check that in the
composition involved in the proof circular components never arise;
then we apply this lemma.

From the Generating Lemma for \SJ\ diagrams it follows that the
homomorphism from \SJo\ to ${(\De,\cirk,I)}$ defined via $\iota$
is onto, and from the Generating Lemma for \SK\ diagrams it
follows that the homomorphism from \SKo\ to ${(\De\times
N,\cirk,(I,0))}$ defined via $\kappa$ is onto. In the next section
we show that these homomorphisms are also one-one.

\section{\SKo\ and \SJo\ are monoids of \SK\ and \SJ\ diagrams}

It is easy to see that the following holds.

\vspace{2ex}

\noindent {\sc Remark.} \emph{Suppose $t$ is the following term of
\SKo\ in normal form:}

\[
c^l\gk{i_p,j_p}\ldots\gk{i_1,j_1}\mj[k_1,l_1]\ldots[k_q,l_q]\dk{m_1,n_1}\ldots\dk{m_r,n_r}.
\]

\noindent \emph{Then}

\nav{(i)}{${\kappa(t)=(\iota(t),l')}$ iff ${l'=l}$;}

\nav{(ii)}{\emph{if ${p>0}$, then ${\{-j,-(i\pl 1)\}}$ is the
maximal cap of ${\iota(t)}$ iff ${j=j_p}$ and ${i=i_p}$; if
${r>0}$, then ${\{m,n\pl 1\}}$ is the maximal cup of ${\iota(t)}$
iff ${m=m_r}$ and ${n=n_r}$;}}

\nav{(iii)}{\emph{if ${p=r=0}$ and ${q>0}$, then ${\{l,-(k\pl
1)\}}$ is the maximal falling thread of ${\iota(t)}$ iff ${l=l_q}$
and ${k=k_q}$.}}

\vspace{2ex}

In a cup ${\{l,k\}}$ of an \SJ\ diagram we call $l$ (which
according to our convention of Section~3 is lesser than $k$) the
\emph{left end point} of this cup. Analogously, in a cap
${\{-l,-k\}}$ we call $-l$ the \emph{left end point} of this cap.
Then under the assumption for $t$ of the Remark we have that

\vspace{1ex}

\nav{(iv)}{\emph{the numbers ${-j'_1,\ldots,-j'_{p'}}$ are the
left end points of the caps of ${\iota(t)}$ iff ${p=p'}$ and
${j'_1,\ldots,j'_p}$ is ${j_1,\ldots,j_p}$; the numbers
${m'_1,\ldots,m'_{q'}}$ are the left end points of the cups of
${\iota(t)}$ iff ${q=q'}$ and ${m'_1,\ldots,m'_q}$ is
${m_1,\ldots,m_q}$.}}

\vspace{1ex}

To prove (iv) it is enough to see that in ${\iota(\gk{i_p,j_p})}$
we have a transversal thread ${\{j_{p-1},-j_{p-1}\}}$, and
analogously with cups. But we will have no use for (iv) below. We
can then prove the following lemma.

\prop{Auxiliary Lemma}{\\[1ex] \indent$(\cal K)$\quad If $t$ and
$u$ are terms of \SKo\ in normal form and ${\kappa(t)=\kappa(u)}$,
then $t$ and $u$ are the same term.
\\[1ex] \indent$(\cal J)$\quad If $t$ and $u$ are terms of \SKo\
in $\cal J$-normal form and ${\iota(t)=\iota(u)}$, then $t$ and
$u$ are the same term.}

\dkz ($\cal K$) Suppose $t$ and $u$ are terms of \SKo\ in normal
form and ${\kappa(t)=\kappa(u)}$. Since ${\kappa(t)=(\iota(t),l)}$
and ${\kappa(u)=(\iota(u),l')}$, we have ${\iota(t)=\iota(u)}$ and
${l=l'}$.

Suppose $t$ is

\[c^l\gk{i_p,j_p}\ldots\gk{i_1,j_1}\mj[k_1,l_1]\ldots[k_q,l_q]\dk{m_1,n_1}\ldots\dk{m_r,n_r}.
\]

\noindent Then the \emph{weight} ${w(t)}$ of $t$ is ${p\pl q\pl
r}$ if ${l=0}$ and ${p\pl q\pl r\pl 1}$ if ${l>0}$. We prove
($\cal K$) by induction on ${w(t)}$.

If ${w(t)=0}$, then $t$ is $\mj$, and by the Remark, $u$ must be
$\mj$ too.

If ${l>0}$, then, by Remark (i), we must have that $u$ is of the
form ${c^lu'}$ where $u'$ is without circles and in normal form;
$t$ is then of the form ${c^lt'}$ for $t'$ in normal form, and
${w(t')=w(t)\mn 1}$. It is easy to see that ${\iota(t)=\iota(t')}$
and ${\iota(u)=\iota(u')}$, so that we have
${\iota(t')=\iota(u')}$. Hence
$\kappa(t')=(\iota(t'),0)=(\iota(u'),0)=\kappa (u')$, and we can
apply the induction hypothesis.

If ${p>0}$, then, by the Remark, we have that $u$ is of the form
${c^l\gk{i_p,j_p}u''}$, where ${c^lu''}$ is in normal form; $t$ is
then of the form ${c^l\gk{i_p,j_p}t''}$, where ${c^lt''}$ is in
normal form, and ${w(c^lt'')=w(t)\mn 1}$.

We have the following equations in \SKo:

\begin{tabbing}
\hspace{7em}\= $[i\mn 1,j]\dk{i\pl 1}\gk{i,j}$ \= = \= $\mj$,\quad
if $j<i$
\\[1ex]
\> \> $\dk{i\pl 1}\gk{i,j}$\' = \> $\mj$,\quad if $j=i$.
\end{tabbing}

\noindent If $v$ stands for ${[i_p\mn 1,j_p]\dk{i_p\pl 1}}$ or
${\dk{i_p\pl 1}}$, depending on whether ${j_p<i_p}$ or
${j_p=i_p}$, then

\begin{tabbing}
\hspace{7em}\= $\kappa(v)\cirk\kappa(t)$ \= = \=
$\kappa(v)\cirk\kappa(u)$
\\[.5ex]
\> \> $\kappa(vt)$\' = \> $\kappa(vu)$
\\[.5ex]
\> \> $\kappa(c^lt'')$\' = \> $\kappa(c^lu'')$,\quad by the
Soundness Lemma,
\end{tabbing}

\noindent and we can apply the induction hypothesis. We proceed
analogously if ${r>0}$, using the following equations of \SKo:

\begin{tabbing}
\hspace{7em}\= $\dk{j,i}\gk{i\pl 1}\sigma_j\ldots\sigma_{i- 1}$ \=
= \= $\mj$,\quad if $j<i$
\\[1ex]
\> \> $\dk{j,i}\gk{i\pl 1}$\' = \> $\mj$,\quad if $j=i$.
\end{tabbing}

If ${p=r=0}$ and ${q>0}$, then, by the Remark, we have that $u$ is
of the form ${u'''[k_q,l_q]}$ for $u'''$ in normal form; $t$ is
then of the form ${t'''[k_q,l_q]}$ for $t'''$ in normal form, and
${w(t''')=w(t)\mn 1}$.

We have the following equations in \SKo:

\[
[k,l]\sigma_l\ldots\sigma_k=\sigma_l\ldots\sigma_k[k,l]=\mj.
\]

\noindent From that we conclude that
${\kappa(t''')=\kappa(u''')}$, and we can apply the induction
hypothesis.

With that the proof of ($\cal K$) is over. The proof of ($\cal J$)
is analogous, with the case ${l>0}$ of the induction step omitted.
\qed

\vspace{2ex}

Finally, we can prove the following.

\prop{Completeness Lemma}{\\[1ex] \indent$(\cal K)$\quad If
${\kappa(t)=\kappa(u)}$, then ${t=u}$ in \SKo.
\\[1ex] \indent$(\cal J)$\quad If ${\iota(t)=\iota(u)}$, then ${t=u}$ in \SJo.}

\dkz ($\cal K$) By the Normal Form Lemma of Section~4, for every
term $t$ and every term $u$ of \SKo\ there are terms $t'$ and $u'$
in normal form such that ${t=t'}$ and ${u=u'}$ in \SKo. By the
Soundness Lemma ($\cal K$) of Section~3, we obtain
${\kappa(t)=\kappa(t')}$ and ${\kappa(u)=\kappa(u')}$, and if
${\kappa(t)=\kappa(u)}$, it follows that
${\kappa(t')=\kappa(u')}$. Then, by the Auxiliary Lemma, $t'$ and
$u'$ are the same term, and hence ${t=u}$ in \SKo. We proceed
analogously for ($\cal J$). \qed

\vspace{2ex}

In the proof of the Auxiliary Lemma we can restrict ourselves to
terms of \SKo\ without cups and caps, and thereby obtain a
Completeness Lemma for the standard presentation of symmetric
groups.

So we have established that the homomorphism from \SKo\ to
${(\De\times N,\cirk,(I,0))}$ defined via $\kappa$ and the
homomorphism from \SJo\ to ${(\De,\cirk,I)}$ defined via $\iota$
are isomorphisms. We can also conclude that for every term $t$ of
\SKo\ there is a \emph{unique} term $t'$ in normal form such that
${t=t'}$ in \SKo. If ${t=t'}$ and ${t=t''}$ in \SKo, then
${t'=t''}$ in \SKo, and hence, by the Soundness Lemma,
${\kappa(t')=\kappa(t'')}$. If $t'$ and $t''$ are in normal form,
we obtain by the Auxiliary Lemma that $t'$ and $t''$ are the same
term. We have an analogous uniqueness for the $\cal J$-normal form
with respect to \SJo.

Let \SKn\ be the submonoid  of \SKo\ made of elements of \SKo\
whose \SK\ diagrams are all the \SK\ diagrams of type ${(n,n)}$.
This monoid corresponds to the multiplicative part of Brauer's
algebra $B_n$ generated out of the basis (see \cite{B37} and
\cite{W88}). How the monoid \SKn\ might be axiomatized (presented
by generators and relations) can be gathered from \cite{Y88} and
\cite{BW89}, where, however, instead of symmetry one finds
braiding.

\section{The maximality of \SKo}

Suppose $t$ and $u$ are terms of \SKo\ such that ${t=u}$ does not
hold in \SKo, and let $\cal X$ be the monoid defined as \SKo\ save
that we assume the additional equation ${t=u}$. We will prove the
following.

\prop{Maximality of \SKo}{For some ${k,l\in \N}$ such that ${k\neq
l}$ we have ${c^k=c^l}$ in $\cal X$.}

\dkz Remember that we have seen in the proof of the Auxiliary
Lemma of the preceding section that for every ${i,j\in \N^+}$ such
that ${j\leq i}$ we have terms $v_1$ and $v_2$ such that
${v_1\gk{i,j}=\mj}$ and ${\dk{j,i}v_2=\mj}$. For terms made
uniquely of crossing terms we also have two-sided inverses. Let
the normal forms of $t$ and $u$ be respectively

\[
\begin{array}{c}
c^l\gk{i_p,j_p}\ldots\gk{i_1,j_1}\mj[k_1,l_1]\ldots[k_q,l_q]\dk{m_1,n_1}\ldots\dk{m_r,n_r},
\\[1ex]
c^{l'}\gk{i'_{p'},j'_{p'}}\ldots\gk{i'_1,j'_1}\mj[k'_1,l'_1]\ldots[k'_{q'},l'_{q'}]
\dk{m'_1,n'_1}\ldots\dk{m'_{r'},n'_{r'}}.
\end{array}
\]

(1) Suppose these normal forms differ only by having ${l\neq l'}$.
Then ${c^l=c^{l'}}$ holds in $\cal X$ by the remark in the
preceding paragraph.

(2) Suppose now that these normal forms differ in some other way,
while either ${l\neq l'}$ or ${l=l'}$. Assume ${l'\leq l}$. Let
$t_1$ be an abbreviation for ${\dk{j_1,i_1}\ldots\dk{j_p,i_p}}$,
while $t_2$ is an abbreviation for

\[
\gk{n_r,m_r}\ldots\gk{n_1,m_1}(\sigma_{l_q}\ldots\sigma_{k_q})
\ldots(\sigma_{l_1}\ldots\sigma_{k_1}).
\]

\noindent It can then be shown  that ${t_1tt_2=c^{l+ p+ r}}$,
while ${t_1ut_2=c^{l'+ k}u'}$, where ${k\leq p\pl r}$ and $u'$ is
in $\cal J$-normal form. If ${k=p\pl r}$, then we do not have that
$u'$ is $\mj$, because otherwise we would be in case (1), as can
be seen from the corresponding \SK\ diagrams.

So in $\cal X$ we have ${c^{l'+ k}u'=c^{l+ p+ r}}$. If $u'$ is
$\mj$, then ${k<p\pl r}$, and so $l'\pl k<l\pl p\pl r$, and we
have proved ${c^{l+ k}=c^{l+ p+ r}}$ in $\cal X$. If, on the other
hand, $u'$ is not $\mj$, then we have the following cases.

(2.1) Let $u'$ be ${\gk{i,j}u''}$. Then in \SKo\ we have

\begin{tabbing}
\hspace{9.5em}\= $\gk{i,j}\dk{j,i}u'$ \= = \=
$\gk{i,j}\dk{j,i}\gk{i,j}u''$
\\[.5ex]
\> \> = \> $cu'$.
\end{tabbing}

\noindent So in $\cal X$ we have

\begin{tabbing}
\hspace{8em}\= $c^{l+ p+ r}\gk{i,j}\dk{j,i}$ \= = \= $c^{l'+
k}\gk{i,j}\dk{j,i}u'$
\\[.5ex]
\> \> = \> $c^{l'+ k+ 1}u'$
\\[.5ex]
\> \> = \> $c^{l+ p+ r+ 1}$,
\end{tabbing}

\noindent and hence we have

\begin{tabbing}
\hspace{6em}\= $c^{l+ p+ r}\dk{i\pl 1}\gk{i,j}\dk{j,i}\gk{i\pl 1}$
\= = \= $c^{l+ p+ r+ 1}\dk{i\pl 1}\gk{i\pl 1}$
\\[.5ex]
\> \> $c^{l+ p+ r}$\' = \> $c^{l+ p+ r+ 2}$.
\end{tabbing}

\noindent We proceed analogously if $u'$ is of the form
${u''\dk{j,i}}$.

(2.2) Let $u'$ be ${\mj u''[k,l]}$. Then in $\cal X$ we have

\[
c^{l+ p+ r+ 1}=c^{l'+ k}\dk{k\pl 1}u'\gk{k\pl 1}.
\]

\noindent Let $u'''$ be the $\cal J$-normal form of ${\dk{k\pl
1}u'\gk{k\pl 1}}$; the term $u'''$ is made uniquely of crossing
terms and $\mj$. For ${k\leq 2n\mn 1}$ we have in $\cal X$

\begin{tabbing}
\hspace{2em}$\dk{1}\dk{3}\ldots\dk{2n\mn 1}c^{l+ p+ r+
1}$\=$\gk{2n\mn 1}\ldots\gk{3}\gk{1}$
\\[.5ex]
\> = \= $\dk{1}\dk{3}\ldots\dk{2n\mn 1}c^{l'+ k}u'''\gk{2n\mn
1}\ldots\gk{3}\gk{1}$
\\[1ex]
\> $c^{l+ p+ r+ 1+ n}$\' = \> $c^{l'+ k+ n'}$
\end{tabbing}

\noindent for ${n'\leq n}$. \qed

\vspace{2ex}

The demonstration in (2.2) could, as a matter of fact, proceed in
different ways. Another such way would be to rely on the following
property of symmetric groups.

Let \So\ be the group generated by crossing terms that satisfies
the monoidal equations and $\sigma$ equations of Section~2. This
is the union of the chain ${\cal S}_1\subset{\cal S}_2\subset{\cal
S}_3\subset\ldots$ of all finite symmetric groups. Then we can
establish the following.

\prop{Maximality of \So}{If ${t=\mj}$ does not hold in \So, and
$\cal X$ is defined as \So\ save that we assume the additional
equation ${t=\mj}$, then for every ${m,n\in\N^+}$ we have
${\sigma_m=\sigma_n}$ in $\cal X$.}

\noindent To prove this we have first the following lemma.

\prop{Lemma 1}{If in $\cal X$ we have ${\sigma_i=\sigma_j}$ for
${i<j}$, then in $\cal X$ we have
${\sigma_i=\sigma_j=\sigma_{j+1}}$.}

\dkz We have
\begin{tabbing}
\hspace{5em}$\sigma_{j+1}$ \= = \=
$\sigma_j\sigma_{j+1}\sigma_j\sigma_{j+1}\sigma_j$,\quad by (1),
($\sigma$2) and ($\sigma$3)
\\[.5ex]
\> = \> $\sigma_i\sigma_{j+1}\sigma_i\sigma_{j+1}\sigma_i$,\quad
since $\sigma_i=\sigma_j$
\\[.5ex]
\> = \> $\sigma_i$,\quad by ($\sigma$), ($\sigma$2) and (1).
\`$\dashv$
\end{tabbing}

\noindent We prove analogously the following lemma.

\prop{Lemma 2}{If in $\cal X$ we have ${\sigma_i=\sigma_j}$ for
${i<j}$, then in $\cal X$ we have
${\sigma_{i-1}=\sigma_i=\sigma_j}$, provided ${i>1}$.}

\noindent Then we prove the Maximality of \So\ as follows.

\vspace{2ex}

\noindent \emph{Proof of the Maximality of \So.} Suppose $t'$ is
the $\cal J$-normal form of $t$, which must be different from
$\mj$. So $t'$ is of the form ${u\sigma_i v}$ where $i$ is the
maximal index of the crossing terms in $t'$; so all indices of
crossing terms in $u$ and $v$ are strictly smaller. From
${t'=\mj}$ it follows that ${\sigma_i=u^{-1}v^{-1}}$, where the
maximal index in ${u^{-1}v^{-1}}$ is strictly smaller than $i$.
Then we have in $\cal X$

\begin{tabbing}
\hspace{11em}$\sigma_{i+1}\sigma_i\sigma_{i+1}$ \= = \=
$\sigma_{i+1}u^{-1}v^{-1}\sigma_{i+1}$
\\[.5ex]
\> $\sigma_i\sigma_{i+1}\sigma_i$\' = \> $u^{-1}v^{-1}$
\\[.5ex]
\> $\sigma_i\sigma_{i+1}\sigma_i$\' = \> $\sigma_i$
\\[.5ex]
\> $\sigma_i$\' = \> $\sigma_{i+1}$.
\end{tabbing}

\noindent It suffices now to apply Lemmata 1 and 2. \qed

\vspace{2ex}

Then we can infer that in the case (2.2) of the proof of the
Maximality of \SKo\ we would have that
${c^{l+p+r}\sigma_1=c^{l'+k}\sigma_2}$ holds in $\cal X$, from
which we infer

\begin{tabbing}
\hspace{11em}$\dk{1}c^{l+p+r}\sigma_1\gk{1}$ \= = \=
$\dk{1}c^{l'+k}\sigma_2\gk{1}$
\\[.5ex]
\> $c^{l+p+r+1}$\' = \> $c^{l'+k}$.
\end{tabbing}

The maximality of \So\ is related to questions concerning the
existence of proper nontrivial subgroups of symmetric groups.

Suppose $t$ and $u$ are terms of \SKo\ such that ${t=u}$ does not
hold in \SJo, and let $\cal X$ be the monoid defined as \SJo\ save
that we assume the additional equation ${t=u}$. Then it is
possible to prove the following.

\prop{Maximality of \SJo}{For every ${k\in\N^+}$ we have
${\gk{k}\dk{k}=\mj}$ in $\cal X$.}

\noindent The proof of that would be an extension of a proof that
may be found in \cite{DP03b}, Section 10, in which we would rely
also on the Maximality of \So. We will not go into this proof,
since we will have no occasion to apply the Maximality of \SJo\ in
this paper.

As a consequence of the Maximality of \SJo\ we have that in $\cal
X$ we have ${\sigma_k=\mj}$ for every ${k\in\N^+}$ (which entails
${\sigma_m=\sigma_n}$, but which need not hold in the $\cal X$ of
the Maximality of \So). In $\cal X$ we also have ${\dk{k}=\dk{k\pl
1}}$ and ${\gk{k}=\gk{k\pl 1}}$, and $\cal X$ is isomorphic to the
group $\Z/n$ for some ${n\in\N}$.

We have considered the matters of this section so as to obtain the
following consequence of the Maximality of \SKo, which we will
apply later. Let $h$ be a monoid homomorphism from \SKo\ to a
monoid $\cal M$. Then \SKo\ is isomorphic to the submonoid
${h(\SKo)}$ of $\cal M$ iff for every ${k,l\in\N}$ such that
${k\neq l}$ we have in $\cal M$ that ${h(c^k)\neq h(c^l)}$.

\section{Symmetric endoadjunctions and self-ad\-junc\-tions}

An {\it adjunction} is a sextuple $\langle{\cal A},{\cal B},
F,G,\varphi,\gamma\rangle$ where, first, $\cal A$ and $\cal B$ are
categories. (Throughout, we deal only with small categories.) This
means that for ${f\!:A\str B}$, ${g\!:B\str C}$ and ${h\!:C\str
D}$ arrows of $\cal A$ we have the equations

\[
\begin{array}{ll}
{\makebox[1cm][l]{$({\mbox{\it cat }} 1)$}} &
{\makebox[6cm][l]{$\mj_B\cirk f=f\cirk\mj_A=f,$}}
\\[1ex]
({\mbox{\it cat }} 2) & h\cirk(g\cirk f)=(h\cirk g)\cirk f,
\end{array}
\]

\noindent and analogously for arrows of ${\cal B}$. Next, $F$, the
\emph{left adjoint}, is a functor from $\cal B$ to $\cal A$, which
means that we have the equations

\[
\begin{array}{ll}
{\makebox[1cm][l]{$({\mbox{\it fun }} 1)$}} &
{\makebox[6cm][l]{$F\mj_B=\mj_{FB},$}}
\\[1ex]
({\mbox{\it fun }} 2) & F(g\cirk f)=Fg\cirk Ff,
\end{array}
\]

\noindent and $G$, the \emph{right adjoint}, is a functor from
$\cal A$ to $\cal B$. Next, $\varphi$, the \emph{counit} of the
adjunction, is a natural transformation from the composite functor
$FG$ from $\cal A$ to $\cal A$ to $Id_{\cal A}$, the identity
functor on $\cal A$, with members (components)
${\varphi_A\!:FGA\str A}$, and $\gamma$, the \emph{unit} of the
adjunction, is a natural transformation from $Id_{\cal B}$ to the
composite functor $GF$ from $\cal B$ to $\cal B$, with members
${\gamma_B\!:B\str GFB}$; this means that we have the equations

\[
\begin{array}{ll}
{\makebox[1cm][l]{$({\mbox{\it nat }} \varphi)$}} &
{\makebox[6cm][l]{$f\cirk\varphi_A=\varphi_B\cirk FGf,$}}
\\[1ex]
({\mbox{\it nat }} \gamma) & GFf\cirk\gamma_A=\gamma_B\cirk f.
\end{array}
\]

\noindent Finally, we have the {\it triangular equations}

\[
\begin{array}{ll}
{\makebox[1cm][l]{$(\varphi\gamma F)$}} &
{\makebox[6cm][l]{$\varphi_{FB}\cirk F\gamma_B=\mj_{FB},$}}
\\[1ex]
(\varphi\gamma G) & G{\varphi_A}\cirk\gamma_{GA}=\mj_{GA}.
\end{array}
\]

\noindent (For the notion of adjunction, and its importance in
category theory, and mathematics in general, see \cite{ML71},
Chapter IV.)

An \emph{endoadjunction} is a quintuple $\langle{\cal A},
F,G,\varphi,\gamma\rangle$ such that $\langle{\cal A},{\cal A},
F,G,\varphi,\gamma\rangle$ is an adjunction.

A \emph{symmetric endoadjunction} is a sextuple $\langle{\cal A},
F,G,\varphi,\gamma, \chi^{GF}\rangle$ such that $\langle{\cal A},
F,G,\varphi,\gamma\rangle$ is an endoadjunction, and with the
definition
\[
\chi_A=_{\df}FF(\varphi_A\cirk\chi^{GF}_A)\cirk
F\chi^{GF}_{FA}\cirk \chi^{GF}_{FFA}\cirk\gamma_{FFA}\!:FFA\str
FFA,
\]
we have that $\chi$ is a natural transformation from $FF$ to $FF$
whose members are self-inverse; this means that we have the
following equations:
\[
\begin{array}{ll}
{\makebox[1cm][l]{(\emph{nat} $\chi$)}} &
{\makebox[6cm][l]{$FFf\cirk\chi_A=\chi_B\cirk FFf,$}}
\\[1ex]
(\chi\chi) & \chi_A\cirk\chi_A=\mj_{FFA}.
\end{array}
\]
We assume, moreover, the equation
\[
\begin{array}{ll}
{\makebox[1cm][l]{($\chi\chi\chi$)}} &
{\makebox[6cm][l]{$\chi_{FA}\cirk
F\chi_A\cirk\chi_{FA}=F\chi_A\cirk\chi_{FA}\cirk F\chi_A,$}}
\end{array}
\]
and, finally, $\chi^{GF}$ is a family of arrows indexed by all the
objects of $\cal A$ such that for every such object $A$ the arrow
${\chi^{GF}_A\!\!:GFA\str FGA}$ is the inverse of
\[
\chi^{FG}_A=_{\df}\;\: GF\varphi_A\cirk
G\chi_{GA}\cirk\gamma_{FGA}\!:FGA\str GFA,
\]
which means that we have the equations
\[
\chi^{GF}_A\cirk\chi^{FG}_A=\mj_{FGA},\quad\quad
\chi^{FG}_A\cirk\chi^{GF}_A=\mj_{GFA}.
\]

Instead of making our assumptions for $\chi$, we could
alternatively make them for the natural transformation $\chi^{GG}$
from $GG$ to $GG$, whose members are defined by
\[
\chi^{GG}_A=_{\df}\;\: GG\varphi_A\cirk
G\chi^{FG}_{GA}\cirk\gamma_{GGA}\!:GGA\str GGA.
\]

When $\langle{\cal A},{\cal B}, F,G,\varphi,\gamma\rangle$ is an
adjunction such that $F$ is left adjoint to $G$ and $\langle{\cal
B},{\cal A}, G,F,\varphi',\gamma'\rangle$ is an adjunction such
that $G$ is left adjoint to $F$, then $\langle{\cal A},{\cal B},
F,G,\varphi,\gamma, \varphi',\gamma'\rangle$ will be called a
\emph{bijunction}. It can then be checked that for every symmetric
endoadjunction $\langle{\cal A},F,G,\varphi,\gamma,
\chi^{GF}\rangle$ we have a bijunction $\langle{\cal A},{\cal
A},F,G,\varphi,\gamma, \varphi',\gamma'\rangle$ where $\varphi'$
and $\gamma'$ are defined by

\[
\begin{array}{rl}
 \varphi'_A\!\!\!\!&=_{\df}\;\:
\varphi_A\cirk\chi^{GF}_A,
\\[1ex]
 \gamma'_A\!\!\!&=_{\df}\;\: \chi^{GF}_A\cirk\gamma_A.
\end{array}
\]

Another alternative definition of symmetric endoadjunction is that
it is a bijunction $\langle{\cal A},{\cal A},F,G,\varphi,\gamma,
\varphi',\gamma'\rangle$ together with a natural transformation
$\chi$ from $FF$ to $FF$ such that the equations ($\chi\chi$) and
($\chi\chi\chi$) are satisfied. With this definition, we define
$\chi^{GF}$ in the following manner:

\[
\chi^{GF}_A=_{\df}\;\:\varphi'_{FGA}\cirk G\chi_{GA}\cirk
GF\gamma'_A.
\]

\noindent This definition may be derived from \cite{T89} (see also
\cite{T94}, Chapter 1).

In every symmetric endoadjunction we have the equations

\[
\begin{array}{ll}
{\makebox[1cm][l]{($\chi\varphi$)}} &
{\makebox[6cm][l]{$\varphi_{FA}\cirk
F\chi^{FG}_A=F\varphi_A\cirk\chi_{GA},$}}
\\[1ex]
(\chi\gamma) & \chi^{FG}_{FA}\cirk
F\gamma_A=G\chi_A\cirk\gamma_{FA}.
\end{array}
\]

A \emph{self-adjunction} is a quadruple $\langle{\cal
A},F,\varphi,\gamma\rangle$ such that $\langle{\cal A},
F,F,\varphi,\gamma\rangle$ is an endoadjunction. (This notion was
considered in \cite{FY89}, Section 4.1, and \cite{DP03a}.) We call
the functor $F$ in a self-adjunction \emph{self-adjoint}. Note
that in a self-adjunction the equation ($\varphi\gamma G$) can be
replaced by

\[
\begin{array}{ll}
{\makebox[1cm][l]{($\varphi\gamma$)}} &
{\makebox[6cm][l]{$\varphi_{FA}\cirk
F\gamma_A=F\varphi_A\cirk\gamma_{FA}.$}}
\end{array}
\]

A \emph{symmetric self-adjunction} is a quintuple $\langle{\cal
A},F,\varphi,\gamma,\chi\rangle$ where $\langle{\cal
A},F,\varphi,\gamma\rangle$ is a self-adjunction, $\chi$ is a
natural transformation from $FF$ to $FF$ such that the equations
($\chi\chi$) and ($\chi\chi\chi$) are satisfied, and where,
moreover, we have

\[
\begin{array}{llll}
{\makebox[1cm][l]{($\chi\varphi 1$)}} &
{\makebox[4,5cm][l]{$\varphi_A\cirk\chi_A=\varphi_A,$}} &
{\makebox[1cm][l]{($\chi\gamma 1$)}} &
{\makebox[4,5cm][l]{$\chi_A\cirk\gamma_A=\gamma_A,$}}
\\[1ex]
{\makebox[1cm][l]{($\chi\varphi 2$)}} &
{\makebox[4,5cm][l]{$\varphi_{FA}\cirk F\chi_A=F\varphi_A\cirk
\chi_{FA},$}} & {\makebox[1cm][l]{($\chi\gamma 2$)}} &
{\makebox[4,5cm][l]{$\chi_{FA}\cirk
F\gamma_A=F\chi_A\cirk\gamma_{FA}.$}}
\end{array}
\]

Alternatively, a symmetric self-adjunction may be defined as a
symmetric endoadjunction where $F$ and $G$ are the same functor,
and where

\[
\chi^{GF}_A=\chi^{FG}_A=\chi_A, \quad\quad \varphi'_A=\varphi_A,
\quad\quad \gamma'_A=\gamma_A.
\]

\noindent The equations ($\chi\varphi 1$) and ($\chi\gamma 1$) are
obtained from the definitions of $\varphi'$ and $\gamma'$, while
the equations ($\chi\varphi 2$) and ($\chi\gamma 2$) are instances
of ($\chi\varphi$) and ($\chi\gamma$) above. Note that
($\chi\varphi 2$) is obtained from ($\varphi\gamma$) by replacing
$\gamma$ by $\chi$, while ($\chi\gamma 2$) is obtained from the
same equation ($\varphi\gamma$) by replacing $\varphi$ by $\chi$.

Let $\kappa_A$ be an abbreviation for
${\varphi_A\cirk\gamma_A\!:A\str A}$. Then in every
self-adjunction, for $f\!:A\str B$ we have that

\[
f\cirk\kappa_A=\kappa_B\cirk f,
\]

\noindent and in every symmetric self-adjunction we have that

\[
F\kappa_A=\kappa_{FA}.
\]

\noindent This is derived as follows:

\begin{tabbing}
\hspace{4em}$F\kappa_A\;$\=$=F\varphi_A\cirk\chi_{FA}\cirk\chi_{FA}\cirk
F\gamma_A$, \quad with ({\mbox{\it fun }} 2) and ($\chi\chi$)
\\*[1ex]
\>$=\varphi_{FA}\cirk F\chi_A\cirk F\chi_A\cirk\gamma_{FA}$, \quad
by ($\chi\varphi 2$) and ($\chi\gamma 2$)
\\[1ex]
\>$=\kappa_{FA}$, \quad with ({\mbox{\it fun }} 2) and
($\chi\chi$).
\end{tabbing}

\noindent So every symmetric self-adjunction is a $\cal
K$-adjunction in the sense of \cite{DP03a} (Section 14).

\section{Free symmetric endoadjunctions and self-ad\-junc\-tions}

The \emph{free symmetric endoadjunction} $\langle{\cal
K'},F,G,\varphi,\gamma,\chi^{GF}\rangle$ generated by a single
object, which we denote by $0$, is defined as follows. The objects
of the category $\cal K'$ of this endoadjunction may be identified
with finite, possibly empty, sequences of occurrences of the
letters $F$ and $G$.

An \emph{arrow-term} of $\cal K'$ will be a word $f$ that has as a
\emph{type} the ordered pair ${(A,B)}$, where $A$ and $B$ are
objects of $\cal K'$. That $f$ is of type $(A,B)$ is expressed by
${f\!:A\str B}$. We define the arrow-terms of $\cal K'$
inductively as follows. We stipulate first for every object $A$ of
$\cal K'$ that ${\mj_A\!:A\str A}$, ${\varphi_A\!:FGA\str A}$,
${\gamma_A\!:A\str GFA}$ and ${\chi^{GF}_A\!:GFA\str FGA}$ are
arrow-terms of $\cal K'$. Next, if ${f\!:A\str B}$ is an
arrow-term of $\cal K'$, then ${Ff\!:FA\str FB}$ and ${Gf\!:GA\str
GB}$ are arrow-terms of $\cal K'$, and if ${f\!:A\str B}$ and
${g\!:B\str C}$ are arrow-terms of $\cal K'$, then ${(g\cirk
f)\!:A\str C}$ is an arrow-term of $\cal K'$. As usual, we do not
write parentheses in ${(g\cirk f)}$ when they are not essential.

On these arrow-terms we impose the equations of symmetric
endoadjunctions. Formally we take the smallest equivalence
relation $\equiv$ on the arrow-terms of $\cal K'$ satisfying,
first, congruence conditions with respect to $F$, $G$ and $\cirk$,
namely,

\[
\begin{array}{l}
{\mbox{\rm if }} f\equiv g,\quad {\mbox{\rm then }} Ff\equiv Fg
{\mbox{\rm{ and }}} Gf\equiv Gg,
\\[1ex]
{\mbox{\rm if }} f_1\equiv f_2\;{\mbox{\rm and }}g_1\equiv g_2,
\quad {\mbox{\rm then }} g_1\cirk f_1\equiv g_2\cirk f_2,
\end{array}
\]

\noindent provided $g_1\cirk f_1$ and $g_2\cirk f_2$ are defined,
and, secondly, the conditions obtained from the equations of
symmetric endoadjunctions by replacing the equality sign $=$ by
$\equiv$. Then we take the equivalence classes of arrow-terms as
arrows, with the obvious sources and targets, all arrow-terms in
the same class having the same type. On these equivalence classes
we define \mj, $\varphi$, $\gamma$, $\chi^{GF}$, $F$, $G$ and
$\cirk$ in the obvious way. This defines the category $\cal K'$,
in which we have clearly a symmetric endoadjunction.

The category $\cal K'$ has the following universal property. If
$s$ maps the object $0$ to an arbitrary object of the category
$\cal A$ of an arbitrary symmetric endoadjunction, then there is a
unique functor $S$ of symmetric endoadjunctions (defined in the
obvious way, so that the structure of symmetric endoadjunctions is
preserved) such that $S$ maps $0$ to $s(0)$. This property
characterizes $\cal K'$ up to isomorphism with a functor of
symmetric endoadjunctions. This justifies calling {\it free} the
symmetric endoadjunction of $\cal K'$.

The \emph{free symmetric self-adjunction} generated by a single
object, whose category we will call $\cal K$, is defined
analogously. The objects of $\cal K$, which are finite, possibly
empty, sequences of occurrences of the letter $F$ preceding $0$,
may be identified with the natural numbers (cf.\ \cite{DP03a},
Section 11).

\section{$\cal K$ and \SKo}

Let $F^0$ be the empty sequence, and let $F^{k+1}$ be $F^kF$. On
the arrows of the category $\cal K$ of the free symmetric
self-adjunction, we define a total binary operation $\ast$ based
on composition of arrows in the following manner. For ${f\!:m\str
n}$ and ${g\!:k\str l}$,
\[
g\ast f=_{\df}\left\{
\begin{array}{ll}
g\cirk F^{k-n}f & {\mbox{\rm if }}n\leq k
\\[1ex]
F^{n-k}g\cirk f & {\mbox{\rm if }}k\leq n.
\end{array}
\right.
\]

Next, let $f\eqk g$ iff there are $k,l\in\N$ such that $F^kf=F^lg$
in $\cal K$. It is easy to check that $\eqk$ is an equivalence
relation on the arrows of $\cal K$, which satisfies moreover

\begin{tabbing}
\hspace{3em}\=$({\mbox{\it
congr}}\;\ast)$\hspace{2em}\=${\mbox{\rm if }}f_1\eqk f_2\;
{\mbox{\rm and }} g_1\eqk g_2, \; {\mbox{\rm then }}g_1\ast
f_1\eqk g_2\ast f_2,$
\\[1ex]
\>$({\mbox{\it congr}}\;F)$\>${\mbox{\rm if }}f\eqk g, \;
{\mbox{\rm then }}Ff\eqk Fg.$
\end{tabbing}

For every arrow $f$ of $\cal K$, let $[f]$ be $\{g\mid f\eqk g\}$,
and let $\cal K^{\ast}$ be $\{[f]\mid f\; {\mbox{\rm is an arrow
of $\cal K$}}\}$. With

\[
\begin{array}{l}
\mj=_{\df}\;[\mj_0],
\\[.1cm]
[g][f]=_{\df}\;[g\ast f],
\end{array}
\]

\noindent we can check that $\cal K^{\ast}$ is a monoid. We will
show that this monoid is isomorphic to the monoid \SKo.

Consider the map $\psi$ from the arrow-terms of $\cal K$ to the
terms of \SKo\ defined inductively by

\begin{tabbing}
\hspace{12em}\=$\xi(\mj)$ \hspace{1.5em}\= is\quad
\=$\mj_0$,\kill

\>$\psi(\mj_n)$\> is \> $\mj$,
\\[.5ex]
\>$\psi(\varphi_n)$\> is \> $\dk{n\pl 1}$,
\\[.5ex]
\>$\psi(\gamma_n)$\> is \> $\gk{n\pl 1}$,
\\[.5ex]
\>$\psi(\chi_n)$\> is \> $\sigma_{n+1}$,
\\[.5ex]
\>$\psi(Ff)$\> is \> $\psi(f)$,
\\[.5ex]
\>$\psi(g\cirk f)$\> is \> $\psi(g)\psi(f)$.
\end{tabbing}

We can easily establish the following by induction on the length
of $f$.

\prop {\sc Remark I}{For every arrow-term ${f\!:m\str n}$ of $\cal
K$, in the $\cal SK$ diagram $\kappa(\psi(f))=(D,k)$ the $\cal SJ$
diagram $D$ is of type $(m,n)$.}

We also have the following.

\prop {\sc Remark II}{If in the $\cal SK$ diagram
$\kappa(t)=(D,k)$ the $\cal SJ$ diagram $D$ is of type $(m,n)$,
then $\kappa(t\dk{m\pl 1})=\kappa(\dk{n\pl 1}t)$,
$\kappa(t\gk{m\pl 1})=\kappa(\gk{n\pl 1}t)$ and
$\kappa(t\sigma_{m+1})=\kappa(\sigma_{n+1}t)$.}

\vspace{2ex}

Then we can prove the following lemma.

\prop {$\psi$ Lemma}{If $f=g$ in $\cal K$, then $\psi(f)=\psi(g)$
in \SKo.}

\dkz We proceed by induction on the length of the derivation of
$f=g$ in $\cal K$. All the cases are quite straightforward except
when $f=g$ is an instance of ({\it nat}~$\varphi$), ({\it
nat}~$\gamma$) or ({\it nat}~$\chi$), where we use Remarks I and
II. In case $f=g$ is an instance of a triangular equation, we use
the equation ({\it cup-cap}) of Section~2; for ($\chi\chi$) we use
($\sigma$2), for ($\chi\chi\chi$) we use ($\sigma$3), for
($\chi\varphi$1) we use ($\sigma$-\emph{cup}~3), for
($\chi\varphi$2) we use ($\sigma$-\emph{cup}~4), for
($\chi\gamma$1) we use ($\sigma$-\emph{cap}~3), for
($\chi\gamma$2) we use ($\sigma$-\emph{cap}~4).
\mbox{\hspace{2em}} \qed

\vspace{2ex}

As an immediate corollary we have that if $f\eqk g$, then
$\psi(f)=\psi(g)$ in \SKo. Hence we have a map from $\cal
K^{\ast}$ to \SKo, which we also call $\psi$, defined by
$\psi([f])=\psi(f)$. Since

\[
\begin{array}{l}
\psi([\mj_0])=\psi(\mj_0)=\mj,
\\[1ex]
\psi([g\ast f])=\psi(g)\psi(f),
\end{array}
\]

\noindent this map is a monoid homomorphism.

Consider next the map $\xi$ from the terms of \SKo\ to the
arrow-terms of $\cal K$ defined inductively by

\begin{tabbing}
\hspace{12em}\=$\xi(\mj)$ \quad\= is\quad  \=$\mj_0$,
\\[.5ex]
\>$\xi(\dk{k})$ \>is\>$\varphi_{k-1}$,
\\[.5ex]
\>$\xi(\gk{k})$ \>is\>$\gamma_{k-1}$,
\\[.5ex]
\>$\xi(\sigma_k)$ \>is\>$\chi_{k-1}$,
\\[.5ex]
\>$\xi(tu)$ \>is\>$\xi(t)\ast\xi(u)$.
\end{tabbing}

\noindent Then we establish the following lemmata.

\prop {$\xi$ Lemma}{If $t=u$ in \SKo, then $\xi(t)\eqk\xi(u)$.}

\dkz We proceed by induction on the length of the derivation of
$t=u$ in \SKo. The cases where $t=u$ is an instance of (1) and (2)
are quite straightforward. For ({\it cup}) or ({\it cup-cap}~1),
by ({\it nat}~$\varphi$) we have

\[
\begin{array}{ll}
\varphi_{k-1}\cirk
F^{k-j+2}\varphi_{j-1}\!\!\!\!\!&=F^{k-j}\varphi_{j-1}\cirk\varphi_{k+1},
\\[1ex]
\varphi_{k+1}\cirk F^{k-j+2}\gamma_{j-1}\!\!\!\!\!&=
F^{k-j}\gamma_{j-1}\cirk \varphi_{k-1}.
\end{array}
\]

\noindent We proceed analogously for ({\it cap}) and ({\it
cap-cup}~1) by using ({\it nat}~$\gamma$). For \mbox{({\it
cup-cap})}, by the triangular equations we have

\[
\begin{array}{ll}
F\varphi_{i-1}\cirk\gamma_i=\mj_i,
\\[1ex]
\varphi_{i-1}\cirk F\gamma_{i-2}=\mj_{i-1}.
\end{array}
\]

\noindent For ($\sigma$), by ({\it nat}~$\chi$) we have

\[
\chi_{k+1}\cirk
F^{k-j+2}\chi_{j-1}=F^{k-j+2}\chi_{j-1}\cirk\chi_{k+1},
\]

\noindent and for ($\sigma$2) and ($\sigma$3) we use ($\chi\chi$)
and ($\chi\chi\chi$) respectively. For ($\sigma$-\emph{cup}~1), by
({\it nat}~$\varphi$) we have

\[
\varphi_{k+1}\cirk
F^{k-j+2}\chi_{j-1}=F^{k-j}\chi_{j-1}\cirk\varphi_{k+1},
\]

\noindent and for ($\sigma$-\emph{cup}~2), by ({\it nat}~$\chi$)
we have

\[
F^{k-j+2}\varphi_{j-1}\cirk\chi_{k+1}=\chi_{k-1}\cirk
F^{k-j+2}\varphi_{j-1}.
\]

\noindent We proceed analogously for ($\sigma$-\emph{cap}~1) and
($\sigma$-\emph{cap}~2) by using ({\it nat}~$\gamma$) and ({\it
nat}~$\chi$) respectively. Finally, for ($\sigma$-\emph{cup}~3) we
use ($\chi\varphi$1), for ($\sigma$-\emph{cup}~4) we use
($\chi\varphi$2), for ($\sigma$-\emph{cap}~3) we use
($\chi\gamma$1), and for ($\sigma$-\emph{cap}~4) we use
($\chi\gamma$2).

We have already established that $\eqk$ is an equivalence relation
that satisfies ({\it congr}~$\ast$) and ({\it congr}~$F$). So the
lemma follows. \qed

\prop {$\xi\psi$ Lemma}{For every arrow-term $f$ of $\cal K$ we
have $\xi(\psi(f))\eqk f$.}

\dkz We proceed by induction on the length of $f$. We have

\begin{tabbing}
\hspace{10em}\=$\xi(\psi(\mj_n))$\quad\=is\hspace{1.5em}\=$\mj_0\;\eqk\;\mj_n$,
\\[.5ex]
\>$\xi(\psi(\varphi_n))$\>is\>$\varphi_n$,
\\[.5ex]
\>$\xi(\psi(\gamma_n))$\>is\>$\gamma_n$,
\\[.5ex]
\>$\xi(\psi(\chi_n))$\>is\>$\chi_n$,
\\[.1cm]
\>$\xi(\psi(Ff))$\>is\>$\xi(\psi(f))$
\\
\>\>                $\eqk$\>$ f$,\quad by the induction hypothesis
\\
\>\>                $\eqk$\> $Ff$,
\\[.1cm]
\>$\xi(\psi(g\cirk f))$\>is\>$ \xi(\psi(g))\ast\xi(\psi(f))$
\\
\>\> $\eqk$\>$            g\cirk f$,
\end{tabbing}

\vspace{-1ex}

\noindent by the induction hypothesis, ({\mbox{\it congr
}}~$\ast$) and the definition of $\ast$. \qed

\vspace{2ex}

By a straightforward induction we can prove also the following
lemma.

\prop {$\psi\xi$ Lemma}{ For every term $t$ of \SKo\ we have that
$\psi(\xi(t))$ is $t$.}

\noindent This establishes that $\cal K^{\ast}$ and \SKo\ are
isomorphic monoids.

Next we prove the following lemma.

\prop {$\cal K$ Cancellation Lemma}{In $\cal K$, if $Ff=Fg$, then
$f=g$.}

\dkz If for ${f,g\!:m\str n}$ we have $m>0$, then

\[
\begin{array}{rl}
\varphi_n\cirk FFf\cirk F\gamma_{m-1} \!\!\!\!& =\, \varphi_n\cirk
FFg\cirk F\gamma_{m-1}
\\[.5ex]
f\cirk\varphi_m\cirk F\gamma_{m-1}\!\!\!\!& =\,
g\cirk\varphi_m\cirk F\gamma_{m-1},\;\;\; {\mbox{\rm by ({\it
nat}~$\varphi$)}}
\\[.5ex]
f \!\!\!\!& =\, g,\;\;\;{\mbox{\rm by ($\varphi\gamma F$)}}.
\end{array}
\]

\noindent If $n>0$, then

\[
\begin{array}{rl}
F\varphi_{n-1}\cirk FFf\cirk\gamma_m \!\!\!\!& =\,
F\varphi_{n-1}\cirk FFg\cirk\gamma_m
\\[.5ex]
F\varphi_{n-1}\cirk\gamma_n\cirk f\!\!\!\!& =\,
F\varphi_{n-1}\cirk\gamma_n\cirk g,\;\;\; {\mbox{\rm by ({\it
nat}~$\gamma$)}}
\\[.5ex]
f \!\!\!\!& =\, g,\;\;\;{\mbox{\rm by ($\varphi\gamma G$)}}.
\end{array}
\]

\noindent If $m=n=0$, then $f$ is equal either to $\mj_0$ or to
$\varphi_0\cirk f'\cirk\gamma_0$, and $g$ is equal either to
$\mj_0$ or to $\varphi_0\cirk g'\cirk\gamma_0$.

If $f=g=\mj_0$, we are done. If $f=\mj_0$ and $g=\varphi_0\cirk
g'\cirk\gamma_0$, then $Ff=Fg$ in $\cal K$ is excluded. From
$Ff=Fg$ in $\cal K$ we obtain
$\kappa(\psi(f))=(I,0)=\kappa(\psi(g))$, which is excluded
because, by Remark~I, in $\kappa(\psi(g'))=(D,k')$ the $\cal SJ$
diagram $D$ must be of type $(2,2)$, from which we obtain that in
$\kappa(\psi(g))=(I,k)$ we have $k\geq 1$. We deal analogously
with the case when $f=\varphi_0\cirk f'\cirk\gamma_0$ and
$g=\mj_0$.

If $f=\varphi_0\cirk f''$ and $g=\varphi_0\cirk g''$, where $f''$
is $f'\cirk\gamma_0$ and $g''$ is $g'\cirk\gamma_0$, then from
$\kappa(\psi(f))=\kappa(\psi(g))$ we conclude
$\kappa(\psi(f''))=\kappa(\psi(g''))$. We can do this because
$\kappa(\psi(f))=\kappa(\psi(g))=(I,k)$ for some $k\geq 1$, while
$\kappa(\psi(f''))=\kappa(\psi(g''))=(\Lambda_1,k\mn 1)$. So
$\psi(f'')=\psi(g'')$, and, by the $\xi\psi$-Lemma, $f''\eqk g''$.
Since $f''$ and $g''$ are both of type $0\str 2$, for some $l\geq
0$ we have $F^l f''=F^l g''$ in $\cal K$, and, by reasoning as at
the beginning of the proof (in the case when $n>0$), we obtain
$f''=g''$ in $\cal K$, and hence $f=g$ in $\cal K$. \qed

\vspace{2ex}

The $\cal K$ Cancellation Lemma implies that $f\eqk g$ for
$f\!:m\str n$ and $g\!:k\str l$ could be defined by $F^{k-n}f=g$
in $\cal K$ when $n\leq k$, and by $f=F^{n-k}g$ in $\cal K$ when
$k\leq n$. So for arbitrary $f,g\!:m\str n$ we have established
that $f=g$ in $\cal K$ iff $f\eqk g$. This establishes that $\cal
K$ is isomorphic to a category of $\cal SK$ diagrams indexed by
all their types.

We define a $\cal J$\emph{-symmetric self-adjunction} as a
symmetric self-adjunction that satisfies moreover
$\varphi_A\cirk\gamma_A=\mj_A$, i.e., $\kappa_A=\mj_A$, according
to the abbreviation introduced at the end of Section~8. In other
words, a $\cal J$-symmetric self-adjunction is a symmetric
self-adjunction that is a $\cal J$-adjunction in the terminology
of \cite{DP03a} (Section 10). Let $\cal J$ be the free $\cal
J$-symmetric self-adjunction generated by a single object, which
is defined analogously to the category $\cal K$, and let $\cal
J^*$ be the monoid defined out of $\cal J$ as $\cal K^*$ is
defined out of $\cal K$. Then we can establish, by imitating what
we had above, that $\cal J^*$ and \SJo\ are isomorphic monoids. We
can also establish that $\cal J$ is isomorphic to a category of
$\cal SJ$ diagrams indexed by all their types.

\section{Subsided categories}

A \emph{monoidal category} $\cal M$ is a category that has a
special object $I$, a bifunctor $\otimes$ on $\cal M$ (i.e.\ a
functor from ${{\cal M}\times{\cal M}}$ to $\cal M$), and natural
isomorphisms whose members (components) are
\begin{tabbing}
\hspace{5em}$a_{A,B,C}$\=: \= $A\otimes(B\otimes C)\str(A\otimes
B)\otimes C$,
\\[.5ex]
\> $\sigma_A\!\!\!$\': \> $I\otimes A\str A$,
\\[.5ex]
\> $\delta_A\!\!\!$\': \> $A\otimes I\str A$,
\end{tabbing}
which satisfy the usual coherence equations (see \cite{ML71},
Section VII.1).

A \emph{symmetric monoidal} category is a monoidal category that
has in addition a natural isomorphism whose members are
\[
s_{A,B}\!:A\otimes B\str B\otimes A,
\]
which satisfy the usual coherence equations (see \cite{ML71},
Section VII.7).

A \emph{symmetric monoidal closed} category is a symmetric
monoidal category $\cal M$ such that for every object $A$ there is
a functor ${[A,\underline{\;\;\;}\,]}$ from $\cal M$ to $\cal M$
right-adjoint to the functor ${A\otimes\underline{\;\;\;}}$, where
for an arrow $f$ of $\cal M$ the arrow ${A\otimes f}$ is
${\mj_A\otimes f}$. The counit of this adjunction has the members

\[
\varepsilon_{A,B}\!:A\otimes[A,B]\str B,
\]

\noindent and the unit has the members

\[
\eta_{A,B}\!:B\str[A,A\otimes B].
\]

\noindent In every symmetric monoidal closed category $\cal M$ we
have a functor ${[\,\underline{\;\;\;}\,,\underline{\;\;\;}\,]}$
from ${{\cal M}^{op}\times{\cal M}}$ to $\cal M$ where for
${f\!:A\str B}$ and ${g\!:C\str D}$

\[
[f,g]=_{df}[A,g]\cirk[A,\varepsilon_{B,C}]\cirk[A,f\otimes\mj_{[B,C]}]\cirk\eta_{A,[B,C]}.
\]

A \emph{compact closed} category is a symmetric monoidal closed
category such that for all objects $A$ and $B$ we have the arrow

\[
\nu_{A,B}\!:[A,B]\str[A,I]\otimes B
\]

\noindent as the inverse of

\[
[\mj_A,\sigma_B\cirk(\varepsilon_{A,I}\otimes\mj_B)\cirk
a_{A,[A,I],B}]\cirk\eta_{A,[A,I]\otimes B}\!:[A,I]\otimes
B\str[A,B].
\]

Alternatively, we could assume the inverse of

\[
[\mj_A,(\varepsilon_{A,C}\otimes\mj_B)\cirk
a_{A,[A,C],B}]\cirk\eta_{A,[A,C]\otimes B}\!:[A,C]\otimes
B\str[A,C\otimes B],
\]

\noindent or of

\begin{tabbing}
\hspace{1em}$[\mj_A,(\mj_B,\varepsilon_{A,C})\cirk
a^{-1}_{B,A,[A,C]}\cirk(s_{A,B}\otimes\mj_{[A,C]})\cirk
a_{A,B,[A,C]}]\cirk\eta_{A,B\otimes[A,C]}\!:$
\\[.5ex]
\` $B\otimes[A,C]\str[A,B\otimes C]$.
\end{tabbing}

\noindent As a matter of fact (as noted in \cite{KL80}, p.~194) it
would be enough to have the arrows $\nu_{A,A}$, because
$\nu_{A,B}$ and all the inverses above can be defined in terms of
$\nu_{A,A}$.

A monoidal category is \emph{strictly monoidal} when for all
objects $A$, $B$ and $C$

\begin{tabbing}
\hspace{11em}$A\otimes(B\otimes C)$ \= = \= $(A\otimes B)\otimes
C$,
\\[.5ex]
\> $I\otimes A$\' = \> $A\otimes I=A$,
\end{tabbing}

\noindent while $a_{A,B,C}=\mj_{A\otimes(B\otimes C)}$ and
${\sigma_A=\delta_A=\mj_A}$.

A \emph{subsided} category is a compact closed category which is
strictly monoidal, in which

\[
[A,B]=[A,I]\otimes B=A\otimes B,
\]

\noindent and where

\begin{tabbing}
\hspace{2em}\= ($\nu$)\quad\quad \= $\nu_{A,B}$ \= = \=
$\mj_{A\otimes B}$,
\\[1ex]
\> ($\varepsilon$1) \> $\varepsilon_{A,B}$\> = \>
$\varepsilon_{A,B}\cirk(s_{A,A}\otimes\mj_B)$,
\\[1ex]
\> ($\varepsilon$2) \> $\varepsilon_{A\otimes
B,C}=\varepsilon_{B,C}\cirk\varepsilon_{A,B\otimes B\otimes
C}\cirk(\mj_A\otimes s_{B,A}\otimes\mj_{B\otimes C})$.
\end{tabbing}

\noindent It can be shown that in every subsided category we have
the equation

\begin{tabbing}
\hspace{2em}\= ($\nu$)\quad\quad \= $\nu_{A,B}$ \= = \=
$\mj_{A\otimes B}$,\kill

\> ($\eta$1) \> $\eta_{A,B}$ \> = \>
$(s_{A,A}\otimes\mj_B)\cirk\eta_{A,B}$,

\end{tabbing}

\noindent dual to ($\varepsilon$1), which could alternatively be
used for defining subsided categories. We also have

\begin{tabbing}
\hspace{2em}\= ($\nu$)\quad\quad \= $\nu_{A,B}$ \= = \=
$\mj_{A\otimes B}$,\kill

\> ($\eta$2) \> $\eta_{A\otimes B,C}=(\mj_A\otimes
s_{A,B}\otimes\mj_{B\otimes C})\cirk\eta_{A,B\otimes B\otimes
C}\cirk\eta_{B,C}$,

\end{tabbing}

\noindent which could be used instead of ($\varepsilon$2).

A compact closed category may be defined as a symmetric monoidal
category $\cal M$ that has a functor $*$ from ${\cal M}^{op}$ to
$\cal M$ and dinatural transformations whose members are

\begin{tabbing}
\hspace{9em}\= $\textbf {\textit {f}}_A$\=$:A\otimes A^*\str I$,
\\[1ex]
\> $\textbf {\textit {g}}_A$\>$:I\str A^*\otimes A$,
\end{tabbing}

\noindent which satisfy the \emph{compact triangular equations}:

\begin{tabbing}
\hspace{5em}\= $\sigma_A\cirk(\textbf {\textit
{f}}_A\otimes\mj_A)\cirk a_{A,A^*,A}\cirk(\mj_A\otimes\textbf
{\textit {g}}_A)\cirk\delta^{-1}_A\;\;$ \= = \= $\mj_A$,
\\[1ex]
\> $\delta_{A^*}\!\cirk(\mj_{A^*}\!\otimes\textbf {\textit
{f}}_A)\cirk a^{-1}_{A^*,A,A^*}\cirk(\textbf {\textit
{g}}_A\otimes\mj_A)\cirk\sigma^{-1}_{A^*}$ \> = \> $\mj_{A^*}$.
\end{tabbing}
That $\textbf {\textit {f}}$ and $\textbf {\textit {g}}$ are
dinatural means that for every arrow $h\!:A\str B$ we have the
equations
\begin{tabbing}
\hspace{9em}\= $\textbf {\textit {f}}_B\cirk(h\otimes
\mj_{B^*})\;$\=$=\textbf {\textit {f}}_A\cirk(\mj_A\otimes h^*)$,
\\[1ex]
\> $(\mj_{A^*}\otimes h)\cirk\textbf {\textit
{g}}_A$\>$=(h^*\otimes \mj_B)\cirk\textbf {\textit {g}}_B$
\end{tabbing}
(see \cite{ML71}, Section IX.4).

With the definition of compact closed category given previously,
$A^*$ is defined as ${[A,I]}$, the arrow $\textbf {\textit {f}}_A$
is defined as $\varepsilon_{A,I}$, and the arrow $\textbf {\textit
{g}}_A$ is defined as
${\nu_{A,A}\cirk[\mj_A,\delta_A]\cirk\eta_{A,I}}$. With the new
definition, we define ${[A,B]}$ as ${A^*\otimes B}$, we define
${[A,f]}$ as ${\mj_{A^*}\otimes f}$, the arrow $\varepsilon_{A,B}$
is defined as ${\sigma_B\cirk(\textbf {\textit
{f}}_A\otimes\mj_B)\cirk a_{A,A^*,B}}$, the arrow $\eta_{A,B}$ is
defined as ${a^{-1}_{A^*,A,B}\cirk(\textbf {\textit
{g}}_A\otimes\mj_B)\cirk\sigma^{-1}_B}$, and $\nu_{A,B}$ is
defined as ${\delta^{-1}_{A^*}\otimes\mj_B}$.

The notion of compact closed category was introduced in
\cite{K72}, and more precisely in \cite{KL80} (see also
\cite{JS93}, Section~7). A subsided category may alternatively be
defined as a symmetric monoidal category which is strictly
monoidal, in which for every object $A$ the functor
${A\otimes\underline{\;\;\;}}$ is self-adjoint, and in this
adjunction we have the equations ($\varepsilon$1) and
($\varepsilon$2). Categories related to subsided categories, but
which instead of symmetry have braiding, were investigated in
\cite{Y88} and \cite{T89} (see also \cite{FY89}, Section 4.1).

\section{$\cal K$  as a subsided category}

To define the free subsided category $\cal S$ generated by a
single object, which we denote by $p$, we rely on the alternative
definition given at the end of the preceding section. The objects
of $\cal S$, which are $I$, $p$, ${p\otimes p}$, ${p\otimes
p\otimes p}$, $\ldots$, may be identified with the natural numbers
$0,1,2,3,\ldots,$ the operation $\otimes$ being addition.

The \emph{arrow-terms} of $\cal S$ are defined inductively as
follows. For all objects $A$ and $B$ of $\cal S$ we have that
${\mj_A\!:A\str A}$, ${s_{A,B}\!:A\otimes B\str B\otimes A}$ (note
that here ${A\otimes B=B\otimes A}$),
${\varepsilon_{A,B}\!:A\otimes A\otimes B\str B}$ and
${\eta_{A,B}\!:B\str A\otimes A\otimes B}$ are arrow-terms of
$\cal S$. Next, if ${f\!: A\str B}$ and ${g\!: C\str D}$ are
arrow-terms of $\cal S$, then $f{\otimes g\!: A\otimes C\str
B\otimes D}$ is an arrow-term of $\cal S$, and if ${f\!: A\str B}$
and ${g\!: B\str C}$ are arrow-terms of $\cal S$, then ${(g\cirk
f)\!: A\str C}$ is an arrow-term of $\cal S$.

On these arrow-terms we impose the equations of subsided
categories; namely, the equations

\begin{tabbing}
\hspace{3em}\=(\emph{bifun}~1)\quad \= \kill

\> (\emph{cat}~1)\> $\mj_B\cirk f=f\cirk\mj_A=f$,
\\[.5ex]
\> (\emph{cat}~2)\> $h\cirk(g\cirk f)=(h\cirk g)\cirk f$,
\\[1ex]
\> (\emph{bifun}~1)\> $\mj_A\otimes\mj_B=\mj_{A\otimes B}$,
\\[.5ex]
\> (\emph{bifun}~2)\> $(g_2\cirk g_1)\otimes(f_2\cirk
f_1)=(g_2\otimes f_2)\cirk(g_1\otimes f_1)$,
\\[1ex]
\> (\emph{nat}~$s$)\> $s_{B,D}\cirk(f\otimes g)=(g\otimes f)\cirk
s_{A,C}$,
\\[.5ex]
\> ($s$2)\> $s_{B,A}\cirk s_{A,B}=\mj_{A\otimes B}$,
\\[.5ex]
\> ($s$3)\> $s_{A\otimes
B,C}=(s_{A,C}\otimes\mj_B)\cirk(\mj_A\otimes s_{B,C})$,
\\[1ex]
\> (\emph{nat}~$\varepsilon$)\>
$f\cirk\varepsilon_{A,B}=\varepsilon_{A,C}\cirk(\mj_{A\otimes
A}\otimes f)$,
\\[.5ex]
\> (\emph{nat}~$\eta$)\> $\eta_{A,C}\cirk f=(\mj_{A\otimes
A}\otimes f)\cirk\eta_{A,B}$,
\\[1ex]
\> ($\varepsilon\eta$)\> $\varepsilon_{A,A\otimes
B}\cirk(\mj_A\otimes\eta_{A,B})=(\mj_A\otimes\varepsilon_{A,B})\cirk\eta_{A,A\otimes
B}=\mj_{A\otimes B}$,
\end{tabbing}

\noindent plus the equations ($\varepsilon$1) and
($\varepsilon$2). Formally, to get the arrows of $\cal S$ we take
equivalence classes as in Section~9.

We will now show that the category $\cal K$ of the free symmetric
self-adjunction generated by a single object (see Section~9) is
isomorphic to $\cal S$. We define first in $\cal K$ the structure
of a subsided category. The objects of $\cal K$ are natural
numbers, the object $I$ is $0$, and $\otimes$ on objects is
addition, which is associative. We define $\otimes$ on the arrows
of $\cal K$ as follows. For ${f\!: m\str n}$

\[
\mj_k\otimes f=_{df} F^k f,
\]

\noindent while ${f\otimes\mj_k}$ is defined inductively as
follows:

\begin{tabbing}
\hspace{5em}\= $\alpha_n\otimes\mj_k$ \= = \= $\alpha_{n+k}$,\quad
for $\alpha$ being $\mj$, $\varphi$, $\gamma$ and $\chi$,
\\[.5ex]
\> $Ff\otimes\mj_k$\> = \> $F(f\otimes\mj_k)$,
\\[.5ex]
\> $(g\cirk f)\otimes\mj_k=(g\otimes\mj_k)\cirk(f\otimes\mj_k)$.
\end{tabbing}

\noindent For ${f\!: m\str n}$ and ${g\!: k\str l}$ we define
${f\otimes g\! : m\pl k\str n\pl l}$ by

\[
f\otimes g=_{df}(\mj_n\otimes g)\cirk(f\otimes\mj_k).
\]

It remains to check that with these definitions $\cal K$ is a
strictly monoidal category, which is easy to do with the help of
\SK\ diagrams. Hence $\cal K$ satisfies (\emph{bifun}~1) and
(\emph{bifun}~2), and ${f\otimes (g\otimes h)=(f\otimes g)\otimes
h}$.

The arrows $s$ are defined inductively as follows:

\begin{tabbing}
\hspace{5em}\= $s_{0,n}=s_{n,0}=\mj_n$,
\\[.5ex]
\> $s_{n+1,1}$ \= = \=
$(s_{n,1}\otimes\mj_1)\cirk(\mj_n\otimes\chi_0)$,
\\[.5ex]
\> $s_{n,m+1}$ \> = \> $(\mj_m\otimes
s_{n,1})\cirk(s_{n,m}\otimes\mj_1)$.
\end{tabbing}

\noindent With this definition $\cal K$ is a symmetric monoidal
category; it satisfies, namely, the equations (\emph{nat}~$s$),
($s$2) and ($s$3) which is easy to check with the help of \SK\
diagrams.

In $\cal K$ the object ${[A,B]}$ is ${A\otimes B}$, while the
arrows $\varepsilon$ and $\eta$ are defined inductively as
follows:

\begin{tabbing}
\hspace{5em}\= $\varepsilon_{0,n}$ \= = \= $\mj_n$,
\\[.5ex]
\> $\varepsilon_{m+1,n}$ \= = \=
$\varphi_n\cirk(\mj_1\otimes\varepsilon_{m,n}\otimes\mj_1)
\cirk(s_{m,1}\otimes\mj_{n+m+1})$,
\\[1.5ex]
\hspace{5em}\= $\varepsilon_{0,n}$ \= = \= $\mj_n$,\kill

\> $\eta_{0,n}$ \> = \> $\mj_n$,
\\[.5ex]
\> $\varepsilon_{m+1,n}$ \= = \=\kill

\> $\eta_{m+1,n}$ \> = \>
$(s_{1,m}\otimes\mj_{n+m+1})\cirk(\mj_1\otimes\eta_{m,n}\otimes\mj_1)\cirk\gamma_n$.
\end{tabbing}

\noindent With these definitions, and with the help of \SK\
diagrams, we check that the equations (\emph{nat}~$\varepsilon$),
(\emph{nat}~$\eta$), ($\varepsilon\eta$), ($\varepsilon$1) and
($\varepsilon$2) are satisfied, which shows that $\cal K$ is a
subsided category. So we have a functor $\Theta$ from the free
subsided category $\cal S$ to $\cal K$ that maps the generating
object $p$ of $\cal S$ into the number $1$, i.e.\ $F0$ (which is
not the generating object of $\cal K$). It is clear that this
functor is a bijection on objects: if objects in both categories
are conceived as natural numbers, then this bijection is identity.
It remains to show that it is a bijection on arrows too.

Now we define in $\cal S$ the structure of a symmetric
self-adjunction. The functor $F$ is
${p\otimes\underline{\;\;\;}}$, the arrow $\varphi_A$ is
$\varepsilon_{p,A}$, the arrow $\gamma_A$ is $\eta_{p,A}$, and the
arrow $\chi_A$ is ${s_{p,p}\otimes\mj_A}$. Since $\cal K$ is free,
we have a functor $\Theta^{-1}$ from $\cal K$ to $\cal S$ that
maps the generating object $0$ into $I$. On objects $\Theta^{-1}$
is identity, as $\Theta$ was. It remains to check that
${\Theta^{-1}(\Theta(f))=f}$ and ${\Theta(\Theta^{-1}(g))=g}$.

It can be verified that the category ${\cal K}'$ of the free
symmetric endoadjunction generated by a single object (see
Section~9) is compact closed. It is not, however, isomorphic to
the free compact closed category generated by a single object; for
this isomorphism to hold, the latter category would need to be
further strictified.

\section{Multiple symmetric self-adjunctions}

Consider a family ${\{\langle{\cal
A},F^i,\varphi^i,\gamma^i\rangle\mid i\in {\cal I}\}}$ of
self-adjunctions, all in the same category $\cal A$, such that for
every ${i,j\in{\cal I}}$ and every object $A$ of $\cal A$ we have
an arrow ${\chi^{i,j}_A\!:F^iF^jA\str F^jF^iA}$, with the
following equations being satisfied:

\begin{tabbing}
\hspace{5em}\= (\emph{nat}~$\chi^i$)\quad \=
$F^jF^if\cirk\chi^{i,j}_A=\chi^{i,j}_B\cirk F^iF^jf$,
\\[1ex]
\> ($\chi^i\chi^i$)\>
$\chi^{j,i}_A\cirk\chi^{i,j}_A=\mj_{F^iF^jA}$,
\\[1ex]
\> ($\chi^i\chi^i\chi^i$)\> $\chi^{j,k}_{F^iA}\cirk
F^j\chi^{i,k}_A\cirk\chi^{i,j}_{F^kA}=F^k\chi^{i,j}_A\cirk\chi^{i,k}_{F^jA}\cirk
F^i\chi^{j,k}_A$,
\\[2ex]
\hspace{.5em}\= ($\chi^i\varphi$1)\quad\=
$\varphi^i_A\cirk\chi^{i,i}_A=\varphi^i_A$,\hspace{7.5em}\=
($\chi^i\gamma$1)\quad\= $\chi^{i,i}_A\cirk\gamma^i_A=\gamma^i_A$,
\\[1ex]
\> ($\chi^i\varphi$2)\> $\varphi^i_{F^jA}\cirk
F^i\chi^{j,i}_A=F^j\varphi^i_A\cirk\chi^{i,j}_{F^iA}$,\>
($\chi^i\gamma$2)\> $\chi^{j,i}_{F^iA}\cirk
F^j\gamma^i_A=F^i\chi^{i,j}_A\cirk\gamma^i_{F^jA}$.
\end{tabbing}

\noindent It is clear that for every $i$ we have that
${\langle{\cal A},F^i,\varphi^i,\gamma^i,\chi^{i,i}\rangle}$ is a
symmetric self-adjunction. We call our family a \emph{multiple
symmetric self-adjunction indexed by} $\cal I$.

Note that every symmetric self-adjunction is a multiple symmetric
self-adjunction indexed by the set of natural numbers $\N$. The
functor $F^0$ is the identity functor, and $F^{n+1}$ is the
composite functor ${F^nF}$. It can be shown that the category of
the multiple symmetric self-adjunction generated by a single
object and the index set $\cal I$ is isomorphic to the free
subsided category generated by the set of objects $\cal I$. We
have treated above this isomorphism in the case where $\cal I$ is
a singleton, and what should be done for the general case is
easily inferred from that.

An arbitrary subsided category $\cal C$ gives rise to a multiple
symmetric self-adjunction indexed by the set of objects of $\cal
C$. For that we proceed as we did in the preceding section for
$\cal S$, putting instead of $p$ an arbitrary object of $\cal C$.

Consider the arrow-terms of the category $\cal M$ of the multiple
symmetric self-adjunction generated by a single object and the
index set $\cal I$. These arrow-terms can be put in a normal form
analogous to the normal form for terms of \SKo. This normal form
looks as follows:

\[
e_{i_1}\ldots e_{i_s}a_p\ldots a_1d_1\ldots d_qb_1\ldots b_r
\]

\noindent where

\nav{(1)}{${i_j\in{\cal I}}$, ${1\leq j\leq s}$, ${i_j\neq i_k}$,
$e_{i_j}$ is
${\varphi^{i_j}_B\cirk\gamma^{i_j}_B\cirk\ldots\cirk\varphi^{i_j}_B
\cirk\gamma^{i_j}_B}$; these arrow-terms correspond to circles,
which are now indexed by members of $\cal I$;}

\nav{(2)}{$a_j$, ${1\leq j\leq p}$, is an arrow-term corresponding
to a block-cap;}

\nav{(3)}{$b_j$, ${1\leq j\leq r}$, is an arrow-term corresponding
to a block-cup;}

\nav{(4)}{and $d_j$, ${1\leq j\leq q}$, is an arrow-term
corresponding to a block-crossing. This normal form is unique up
to the order of the indices ${i_1\ldots i_s}$.}

The category $\cal M$ can be shown isomorphic to a category of
diagrams analogous to \SK\ diagrams indexed by types, but whose
types are now not ${(n,m)}$ with ${n,m\in\N^+}$, but instead
${(A,B)}$ where $A$ and $B$ are finite sequences of elements of
$\cal I$. These new diagrams have circular components indexed by
elements of $\cal I$. The proof of this isomorphism is analogous
to the proof of isomorphism of the category $\cal K$ of the free
symmetric self-adjunction with a category of \SK\ diagrams indexed
by their types (see Section 10).

We can establish that $\cal M$ is maximal in the following sense.
Suppose that $t$ and $u$ are arrow-terms of $\cal M$ such that
${t=u}$ does not hold in $\cal M$, and let $\cal X$ be the
category defined as $\cal M$ save that we assume the additional
equation ${t=u}$. We can then prove the following.

\prop{Maximality of $\cal M$}{For two arrow-terms $v$ and $v'$ of
$\cal M$ in normal form such that $v$ is \[e_{i_1}\ldots
e_{i_s}a_p\ldots a_1d_1\ldots d_qb_1\ldots b_r,\] $v'$ is
\[e_{i'_1}\ldots e_{i'_{s'}}a'_{p'}\ldots a'_1d'_1\ldots d'_{q'}b'_1
\ldots b'_{r'},\] and ${\{e_{i_j}\mid 1\leq j\leq
s\}\neq\{e_{i'_j}\mid 1\leq j\leq s'\}}$, we have that ${v=v'}$
holds in $\cal X$.}

\dkz If the normal forms of $t$ and $u$ differ in their $e$-parts
as above, then we are done.

Suppose the normal forms of $t$ and $u$ do not differ in their
$e$-parts, but they differ either in their $a$-parts or in their
$b$-parts, and suppose the number of $a$-terms and $b$-terms in
the normal form of $t$ is greater than or equal to the
corresponding number in the normal form of $u$. Then with the
arrow-terms $a$ and $b$, whose diagrams are mirror images of the
diagrams corresponding respectively to the $b$ and $a$ parts of
the normal form of $t$, we have

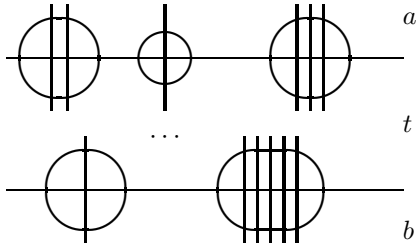
\begin{center}
\begin{picture}(150,110)
{\linethickness{0.05pt} \put(0,30){\line(1,0){150}}
\put(0,80){\line(1,0){150}}}

\put(150,15){\makebox(0,0)[l]{$b$}}
\put(150,55){\makebox(0,0)[l]{$t$}}
\put(150,95){\makebox(0,0)[l]{$a$}}

\thicklines \put(30,10){\line(0,1){40}}
\put(90,10){\line(0,1){40}} \put(95,10){\line(0,1){40}}
\put(100,10){\line(0,1){40}} \put(105,10){\line(0,1){40}}
\put(110,10){\line(0,1){40}}

\put(17,60){\line(0,1){40}} \put(23,60){\line(0,1){40}}
\put(60,60){\line(0,1){40}} \put(110,60){\line(0,1){40}}
\put(115,60){\line(0,1){40}} \put(120,60){\line(0,1){40}}

\put(60,50){\makebox(0,0){$\cdots$}}

\put(30,30){\oval(30,30)} \put(100,30){\oval(40,30)}
\put(20,80){\oval(30,30)} \put(60,80){\oval(20,20)}
\put(115,80){\oval(30,30)}

\end{picture}
\end{center}

\noindent This means that the normal forms of ${b\cirk t\cirk a}$
and ${b\cirk u\cirk a}$ differ in their $e$-parts, because in the
diagram of ${b\cirk u\cirk a}$ we must have less circular
components.

If the normal forms of $t$ and $u$ differ only in their $d$-parts,
then in the diagram corresponding to $t$ we have a transversal
thread

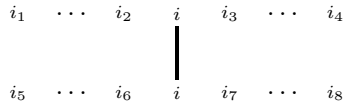
\begin{center}
\begin{picture}(150,60)

\put(0,10){\makebox(0,0){\scriptsize $i_5$}}
\put(40,10){\makebox(0,0){\scriptsize $i_6$}}
\put(60,10){\makebox(0,0){\scriptsize $i$}}
\put(80,10){\makebox(0,0){\scriptsize $i_7$}}
\put(120,10){\makebox(0,0){\scriptsize $i_8$}}
\put(20,10){\makebox(0,0){$\ldots$}}
\put(100,10){\makebox(0,0){$\ldots$}}

\put(0,40){\makebox(0,0){\scriptsize $i_1$}}
\put(40,40){\makebox(0,0){\scriptsize $i_2$}}
\put(60,40){\makebox(0,0){\scriptsize $i$}}
\put(80,40){\makebox(0,0){\scriptsize $i_3$}}
\put(120,40){\makebox(0,0){\scriptsize $i_4$}}
\put(20,40){\makebox(0,0){$\ldots$}}
\put(100,40){\makebox(0,0){$\ldots$}}

\thicklines \put(60,15){\line(0,1){20}}

\end{picture}
\end{center}

\noindent which is missing in the diagram corresponding to $u$,
but in the latter diagram out of the $i$'s above and below go
transversal threads, since the $a$-parts and $b$-parts of the two
diagrams do not differ. Then we have

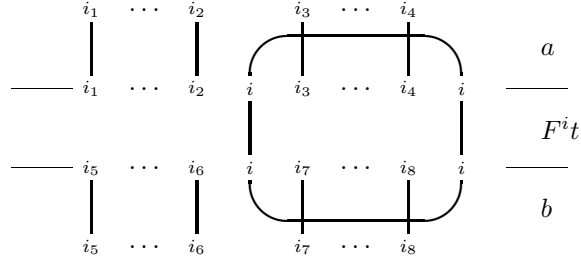
\begin{center}
\begin{picture}(150,110)

{\linethickness{0.05pt} \put(-30,40){\line(1,0){23}}
\put(-30,70){\line(1,0){23}} \put(157,40){\line(1,0){23}}
\put(157,70){\line(1,0){23}}}

\put(0,10){\makebox(0,0){\scriptsize $i_5$}}
\put(40,10){\makebox(0,0){\scriptsize $i_6$}}
\put(80,10){\makebox(0,0){\scriptsize $i_7$}}
\put(120,10){\makebox(0,0){\scriptsize $i_8$}}
\put(20,10){\makebox(0,0){$\ldots$}}
\put(100,10){\makebox(0,0){$\ldots$}}

\put(0,40){\makebox(0,0){\scriptsize $i_5$}}
\put(40,40){\makebox(0,0){\scriptsize $i_6$}}
\put(60,40){\makebox(0,0){\scriptsize $i$}}
\put(80,40){\makebox(0,0){\scriptsize $i_7$}}
\put(120,40){\makebox(0,0){\scriptsize $i_8$}}
\put(140,40){\makebox(0,0){\scriptsize $i$}}
\put(20,40){\makebox(0,0){$\ldots$}}
\put(100,40){\makebox(0,0){$\ldots$}}

\put(0,70){\makebox(0,0){\scriptsize $i_1$}}
\put(40,70){\makebox(0,0){\scriptsize $i_2$}}
\put(60,70){\makebox(0,0){\scriptsize $i$}}
\put(80,70){\makebox(0,0){\scriptsize $i_3$}}
\put(120,70){\makebox(0,0){\scriptsize $i_4$}}
\put(140,70){\makebox(0,0){\scriptsize $i$}}
\put(20,70){\makebox(0,0){$\ldots$}}
\put(100,70){\makebox(0,0){$\ldots$}}

\put(0,100){\makebox(0,0){\scriptsize $i_1$}}
\put(40,100){\makebox(0,0){\scriptsize $i_2$}}
\put(80,100){\makebox(0,0){\scriptsize $i_3$}}
\put(120,100){\makebox(0,0){\scriptsize $i_4$}}
\put(20,100){\makebox(0,0){$\ldots$}}
\put(100,100){\makebox(0,0){$\ldots$}}

\put(170,25){\makebox(0,0)[l]{$b$}}
\put(170,55){\makebox(0,0)[l]{$F^it$}}
\put(170,85){\makebox(0,0)[l]{$a$}}

\thicklines \put(60,45){\line(0,1){20}}
\put(140,45){\line(0,1){20}} \put(0,15){\line(0,1){20}}
\put(40,15){\line(0,1){20}} \put(80,15){\line(0,1){20}}
\put(120,15){\line(0,1){20}}

\put(0,75){\line(0,1){20}} \put(40,75){\line(0,1){20}}
\put(80,75){\line(0,1){20}} \put(120,75){\line(0,1){20}}

\put(100,35){\oval(80,30)[b]} \put(100,75){\oval(80,30)[t]}

\end{picture}
\end{center}

\noindent This means that the normal forms of ${b\cirk F_it\cirk
a}$ and ${b\cirk F_iu\cirk a}$ differ in the $e$-parts, since the
diagram of ${b\cirk F_it\cirk a}$ has an additional circular
component. \qed

\vspace{2ex}

In this section we have not been very precise, because the matters
considered are just complications of matters considered
previously. More precision would involve a considerable number of
these complications, while the ideas would not be essentially new.

\section{The category \Mat}

Let \Mat\ be the skeleton of the category of finite-dimensional
vector spaces over a number field $\cal F$ with linear
transformations as arrows. (A skeleton of a category $\cal C$ is
any full subcategory $\cal C'$ of $\cal C$ such that each object
of $\cal C$ is isomorphic in $\cal C$ to exactly one object of
$\cal C'$; any two skeletons of $\cal C$ are isomorphic
categories, so that, up to isomorphism, we may speak of \emph{the}
skeleton of $\cal C$.)

More precisely, the objects of the category \Mat\ are natural
numbers (the dimensions of our vector spaces), an arrow
${M\!\!:m\str n}$ is an $n\times m$ matrix, composition of arrows
$\cirk$ is matrix multiplication, and the identity arrow
$\mj_n\!:n\str n$ is the identity matrix (i.e.\ the $n\times n$
matrix with 1 on the diagonal and 0 elsewhere). If either $m$ or
$n$ is $0$, then there is only one matrix ${M\!\!:m\str n}$,
called the empty matrix, which is the empty map $\emptyset$ from
$\emptyset$ to $\cal F$ indexed by ${(m,n)}$.

Consider the functor $\otimes$ from $\Mat\times\Mat$ to \Mat\ that
is product on objects and that on arrows, i.e.\ matrices, is the
Kronecker product (see \cite{J53}, Chapter VII.5, pp.\ 211-213).
Let the special object $I$ in \Mat\ be the number $1$. With
$\otimes$ and $I$ the category \Mat\ is a strictly monoidal
category.

Let $S_{m,n}$ be the $nm\times mn$ matrix that for $1\leq i\leq n$
and $1\leq j\leq m$ has the entries

\[
S_{m,n}((i\mn 1)m\pl j, k)=\delta(k,(j\mn 1)n\pl i),
\]

\noindent where $\delta$ is the Kronecker delta. If either $m$ or
$n$ is $0$, then $S_{m,n}$ is the empty matrix $\mj_0$. For
example, $S_{3,2}$ is the matrix

\[
\left[
\begin{array}{cccccc}
  1 & 0 & 0 & 0 & 0 & 0 \\
  0 & 0 & 1 & 0 & 0 & 0 \\
  0 & 0 & 0 & 0 & 1 & 0 \\
  0 & 1 & 0 & 0 & 0 & 0 \\
  0 & 0 & 0 & 1 & 0 & 0 \\
  0 & 0 & 0 & 0 & 0 & 1 \\
\end{array}
\right]
\]

\noindent Then we can check that with $\otimes$, $I$ and $S_{m,n}$
so defined \Mat\ is a symmetric monoidal category.

Let $E_{m,1}$ be the $1\times m^2$ matrix that for $1\leq i,j\leq
m$ has the entries

\[E_{m,1}(1,(i\mn 1)m\pl j)=\delta(i,j).\]

\noindent We define $E_{m,n}$ as the $n\times m^2n$ matrix
$E_{m,1}\otimes \mj_n$. For example, $E_{2,2}$ is

\[
\left[
\begin{array}{cccccccc}
  1 & 0 & 0 & 0 & 0 & 0 & 1 & 0 \\
  0 & 1 & 0 & 0 & 0 & 0 & 0 & 1 \\
\end{array}
\right]
\]

\noindent If ${n=0}$, then $E_{m,n}$ is $\mj_0$, and if ${m=0}$,
then $E_{m,n}$ is the empty matrix ${M\!\!:0\str n}$.

Let $H_{m,n}$ be the transpose of $E_{m,n}$. Taking that $E_{m,n}$
is $\varepsilon_{m,n}$ and $H_{m,n}$ is $\eta_{m,n}$, we can check
that for every $m$ the functor $m\otimes\underline{\;\;\;}$ (where
$m\otimes\underline{\;\;\;}$ on arrows is
$\mj_m\otimes\underline{\;\;\;}$) is self-adjoint. Moreover, the
equations ($\varepsilon$1) and ($\varepsilon$2) are satisfied. So
we can conclude that \Mat\ is a subsided category.

As a matter of fact, we can show that in \Mat\ we have as a
subcategory an isomorphic copy of the skeleton ${\cal M}_s$ of the
category $\cal M$ of the multiple symmetric self-adjunction freely
generated by a single object $0$ and the index set ${\bf N}^+$.
The objects of ${\cal M}_s$ may be identified with multisets of
positive natural numbers.

We have a functor $\textbf {\textit {F}}'$ from ${\cal M}_s$ to
$\cal M$ (which are equivalent categories), and we have a functor
$\textbf {\textit {F}}''$ from $\cal M$ to \Mat\ such that
$\textbf {\textit {F}}''(0)=1$, and, for $p_i$ being the $i$-th
prime number, $\textbf {\textit {F}}''(F^iA)=p^i\textbf {\textit
{F}}''(A)$. A functor \textbf {\textit {F}} from ${\cal M}_s$ to
\Mat\ is obtained as the composite functor ${\textbf {\textit
{F}}''\!\textbf {\textit {F}}'}$. The functor \textbf {\textit
{F}} establishes a one-to-one correspondence between the objects
of ${\cal M}_s$ and the objects of \Mat\ with $0$ omitted.

That the functor \textbf {\textit {F}} is one-one on arrows is
shown as follows. We verify first that $\textbf {\textit {F}}''$
is a faithful functor. This is a consequence of the Maximality of
$\cal M$ of the preceding section, as we will now show.

Suppose $t=u$ does not hold in $\cal M$. If we had $\textbf
{\textit {F}}''(t)=\textbf {\textit {F}}''(u)$ in \Mat, then by
the Maximality of $\cal M$ we would have in \Mat\ that $\textbf
{\textit {F}}''(v)=\textbf {\textit {F}}''(v')$ for $v$ and $v'$
as in the statement of this Maximality. The matrices

\[\textbf {\textit {F}}''(a_p\ldots a_1d_1\ldots d_qb_1\ldots
b_r) \mbox{\rm { and }} \textbf {\textit {F}}''(a_{p'}'\ldots
a_1'd_1'\ldots d_{q'}'b_1'\ldots b_{r'}')\]

\noindent are 0-1 matrices, because $\textbf {\textit {F}}''(a_j)$
is a matrix with at most one $1$ in each row, $\textbf {\textit
{F}}''(b_j)$ is a matrix with at most one 1 in each column, and
$\textbf {\textit {F}}''(d_j)$ is a matrix with at most one 1 in
each row and each column. These matrices are not 0 matrices,
because they all have 1 in the $(1,1)$ entry. The matrix $\textbf
{\textit {F}}''(e_{i_1}\ldots e_{i_s})$ is different from $\textbf
{\textit {F}}''(e_{i_1'}\ldots e_{i_s'})$, because for $p_{i+1}$
the ${i\pl 1}$-th prime number $\textbf {\textit
{F}}''(\varphi^i_B\cirk\gamma^i_B)=p_{i+1}\mj_{\footnotesize\textbf
{\textit {F}}''(B)}$. So $\textbf {\textit {F}}''(t)=\textbf
{\textit {F}}''(u)$ does not hold in \Mat, and hence $\textbf
{\textit {F}}''$ is faithful.

Since $\textbf {\textit {F}}'$ is faithful by definition, we have
that the functor \textbf {\textit {F}}, which is ${\textbf
{\textit {F}}''\!\textbf {\textit {F}}'}$, is faithful, and since
it is one-one on objects, it is one-one on arrows.

That in \Mat\ we have as a subcategory an isomorphic copy of the
category $\cal K$ of the symmetric self-adjunction freely
generated by a single object is a consequence of the fact that in
\Mat\ we have as a subcategory an isomorphic copy of $\cal M$, and
of the fact that in $\cal M$ we have denumerably many isomorphic
copies of $\cal K$. However, this can also be shown directly, and
more easily, by relying on the Maximality of \SKo\ of Section~7.
An alternative proof can be based on \cite{DP03c} and
\cite{DP03d}.

We just said that in $\cal M$ we have denumerably many isomorphic
copies of $\cal K$; namely, there are denumerably many embeddings
$\textbf {\textit {E}}_i$ of $\cal K$ in $\cal M$, such that
$\textbf {\textit {E}}_i(FA)=F^i\textbf {\textit {E}}_i(A)$.
Conversely, there is also a functor \textbf {\textit {E}} from
$\cal M$ to $\cal K$ such that $\textbf {\textit
{E}}(F^iA)=F\textbf {\textit {E}}(A)$. This functor is not
faithful because it maps $\varphi^i_A\cirk\gamma^i_A$ and
$\varphi^j_A\cirk\gamma^j_A$ for $i\neq j$ into
$\varphi_{\footnotesize\textbf {\textit
{E}}(A)}\cirk\gamma_{\footnotesize\textbf {\textit {E}}(A)}$. If,
however, the diagrams corresponding to $f$ and $g$ have no
circular components, and ${\textbf {\textit {E}}(f)=\textbf
{\textit {E}}(g)}$ in $\cal K$, then ${f=g}$ in $\cal M$.

A category analogous to \Mat\ in which we would have an isomorphic
copy of the category $\cal J$ would have as objects natural
numbers and as arrows 0-1 matrices, whose composition is defined
like composition of relations. This is like matrix multiplication
together with $1\pl 1=1$.

\vspace{2ex}

\noindent {\small {\it Acknowledgement$\,$}. Work on this paper
was supported by the Ministry of Science of Serbia (Grant
144013).}

\end{document}